\documentclass[11pt, a4paper]{article}
\usepackage{amsmath}
\usepackage{amsfonts}
\usepackage{bm}
\newcommand\bkE{{\mathbb {E}}}

\newcommand\R{{\mathbb {R}}}
\newcommand\N{{\mathbb {N}}}

\newtheorem{Theorem}{Theorem}
\newtheorem{corollary}[Theorem]{Corollary}

\newtheorem{Lemma}[Theorem]{Lemma}

\newtheorem{Proposition}[Theorem]{Proposition}
\newtheorem{remark}[Theorem]{Remark}

\renewcommand{\P}{\ensuremath{\mathbb {P}}}
\newcommand{\E}{\mathbb{E}}
\newcommand{\LL}{ { \mathbb L}}
\newcommand{\p}{ { \mathbb P} }
\newcommand{\Z}{ { \mathbb Z} }

\newcommand{\F}{\mathcal{F}}

\newcommand\beq{\begin{equation}}
\newcommand\eeq{\end{equation}}

\def\cov{\mathop{\rm Cov}\limits}

\begin{document}

\title{On the Koml\'os, Major and Tusn\'ady  strong approximation for some classes of random iterates}

\author{Christophe Cuny{\footnote{Universit\'e de la Nouvelle-Cal\'edonie, Equipe ERIM. Email: christophe.cuny@univ-nc.nc}}, J\'er\^ome Dedecker{\footnote{Universit\'e Paris Descartes, Sorbonne Paris Cit\'e,  Laboratoire MAP5 (UMR 8145). Email: jerome.dedecker@parisdescartes.fr}} and Florence Merlev\`ede{\footnote{Universit\'{e} Paris-Est, LAMA (UMR 8050), UPEM, CNRS, UPEC. Email: florence.merlevede@u-pem.fr}}}

\maketitle

{\abstract{The famous results of Koml\'os, Major and Tusn\'ady (see \cite{KMT}  
and \cite{Ma76})  state that it is possible to approximate almost  
surely the partial sums of size $n$ of i.i.d. centered random  
variables in ${\mathbb L}^p$ ($p >2$) by a Wiener process with an  
error term of order $o(n^{1/p})$. Very recently, Berkes, Liu and Wu  
\cite{BLW14} extended this famous result to  partial sums associated  
with functions of an i.i.d. sequence, provided  a condition on a  
functional dependence measure in ${\mathbb L}_p$ is satisfied. In this
paper, we adapt the method of Berkes, Liu and Wu to  partial sums of  
functions of random iterates.  Taking  
advantage of the Markovian setting, we shall give new dependent  
conditions, expressed in terms of a natural coupling (in ${\mathbb  
L}^\infty$ or in ${\mathbb L}^1$),
under which the strong approximation result holds with rate  
$o(n^{1/p})$. As we shall see our conditions are well adapted to a  
large variety of models, including left random
walks on $GL_d({\mathbb R})$, contracting iterated random functions,  
autoregressive Lipschitz processes,  and some ergodic Markov chains.
We also provide some examples showing that our ${\mathbb  
L}^1$-coupling condition is in some sense optimal.}}

\section{Introduction}\label{Sec1}

In this paper we shall adapt the approach of Berkes-Liu-Wu 
\cite{BLW14} to certain classes of Markov chains.
To motivate this work, let us describe in detail  the example of the left random walk on $GL_d(\mathbb R)$,  $d \geq 2$ (the group of invertible $d$-dimensional real matrices). 

Let $(\varepsilon_n)_{n \geq 1}$ be independent random matrices taking values in $G= GL_d(\mathbb R)$,  with common distribution $\mu$.  Let $\Vert \cdot \Vert$ be the euclidean norm on ${\mathbb R}^d$. We shall say that $\mu $ has a moment of order $p \geq 1$ if
\begin{equation}\label{Mpmu}
\int_G (\log N(g) )^p \mu(dg) < \infty \, , 
\end{equation}
where $N(g) := \max ( \Vert g \Vert , \Vert g^{-1} \Vert)$. 

Let $A_0={\rm Id}$ and for every $n \geq 1$, $A_n = \varepsilon_n \cdots \varepsilon_1$. 
Recall that if $\mu$ admits a moment of order $1$ then 
\begin{equation}\label{LLN}
\lim_{n \rightarrow \infty} \frac{1}{n} \log \Vert A_n \Vert = \lambda_{\mu} \, \text{ ${\mathbb P}$-a.s.}, 
\end{equation}
where $\lambda_\mu:= \lim_{n\to +\infty} n^{-1}
\E(\log \|\varepsilon_n \cdots \varepsilon_1\|)$ is the so-called first Lyapunov exponent (see for instance \cite{FK}). 
For any $x \in S^{d-1}$, we  want to describe as precisely as possible the asymptotic behavior of the quantity
\begin{equation} \label{lognorm}
\log \Vert  A_n x \Vert  \, .
\end{equation}

The left random walk of law $\mu$ started at $x \in S^{d-1}$ 
is the Markov chain  defined by 
$W_{0,x}:=x$ and $W_{n,x} =\varepsilon_n W_{n-1, x}$ for $n \geq 1$.
As usual, to handle the quantity \eqref{lognorm},  we consider the partial sums associated with the
random variables  $(X_{n,x})_{n \geq 1}$ given by 
\[
X_{n,x} := h (\varepsilon_n, W_{n-1, x} ) \, , \, n \geq 1 \, , 
\]
where for every $g \in G$ and every $y \in {\mathbb R}^{d}-\{0\}$,
\[
 h ( g , y ) = \log \left ( \frac{\Vert g \cdot y \Vert }{ \Vert y \Vert }\right ) \, .
\]
By definition of $h$, and since  $X_{n,x} =h (\varepsilon_n, A_{n-1}x)$, we easily see that, for any $x \in S^{d-1}$, 
\begin{equation}\label{mainidentity}
S_{n,x} = \sum_{k=1}^n X_{k,x} = \log \Vert  A_n x \Vert \, .
\end{equation}
Hence, the asymptotic behavior  of \eqref{lognorm} can 
be deduced from the asymptotic  behavior of partial sums of functions of the Markov chain $W_{n,x}$. 

This problem can be tackled under some assumptions on 
$\mu$ (strong irreducibility and proximality, see subsection 
\ref{LRW} 
for more details) which implies that the chain $(W_n)_{n \geq 0}$ admits an unique invariant measure  $\nu$ defined on the projective 
space $X:=P_{d-1} ({\mathbb R}^d)$  of ${\mathbb R}^d -\{0\}$.
Under these assumptions on $\mu$,  and assuming moreover that $\mu$ has a moment of order $p \in (2, 4)$,  Cuny-Dedecker-Jan \cite{CDJ} proved the following strong approximation result: there exists $\sigma^2\geq 0$ such that,
 for every (fixed) $x \in S^{d-1}$, one can redefine $(
 \log \Vert  A_n x \Vert)_{n \geq 1}$ without changing its distribution on a (richer) probability space on which 
there exist iid random variables $(N_i)_{i \geq 1}$ with common distribution ${\mathcal N} (0, \sigma^2)$, such that,
\begin{equation}\label{strong}
 \log \Vert  A_n x \Vert - n \lambda_{\mu} - \sum_{i=1}^n N_i = o(n^{1/p}\sqrt{\log n}) \, \text{ a.s.}
\end{equation}
If $\mu$ has a moment of order $p=4$,  the same authors showed that this strong approximation holds with a rate of order
$O(n^{1/4} \sqrt{\log (n)} (\log \log n )^{1/4})$. 

To prove \eqref{strong}, Cuny-Dedecker-Jan used a 
martingale approximation (as described for instance 
in Cuny-Merlev\`ede \cite{CM15}), together with some
appropriate upper bounds on the quantities
\begin{equation}\label{sup}
\sup_{\|x\|=1,  \|y\|=1} 
\E \left (|X_{k, x}-X_{k,  y}|\right) \, .
\end{equation}
The main drawback of this approach is that it cannot give 
a better rate than $n^{1/4}$, because it is based on the 
Skorokhod representation theorem for martingales.

On another hand, since the stationary Markov chain $W_n$ is a function 
of the starting point $W_0$ and of the ``innovations" $\varepsilon_1, \cdots, \varepsilon_n$, one can also apply 
the approximation results by Berkes-Liu-Wu (in fact, this is
not completely immediate because it does not fit exactly into the framework described by these authors, and some extra work is required there). Doing so, one can reach a rate of order 
$n^{1/p}$ for any $p>2$, but only by assuming that $\mu$
has a moment of order $q(p)>p$. More precisely, their functional measure of dependence in ${\mathbb L}_p$, say $\delta_{k,p}$,  can be bounded by $\sup_{\|x\|=1,  \|y\|=1} 
\Vert X_{k, x}-X_{k,  y} \Vert_p$. Hence, applying Proposition 3 in \cite{CDJ}, one can see that condition (2.3) in \cite{BLW14} is satisfied provided  $\mu$ has at least a moment of order $(5p/2) -1$.  This is somewhat surprising:
on the one hand, one can go beyond the rate of order $n^{1/4}$, and on the other hand we need stronger assumptions than in  Cuny-Dedecker-Jan  \cite{CDJ} to get the
rate $n^{1/p}$ when $p \in (2,4)$. 

This gave us a strong motivation to understand completely 
the proof by Berkes-Liu-Wu \cite{BLW14}, and to see whether it 
is possible to take advantage of the Markovian setting to get 
the rate $n^{1/p}$ in \eqref{strong} under a moment of order $p$, for any $p>2$.
As we shall see in this paper, the answer is positive.

As already mentioned, in the case of the left random walk on $GL_d(\R)$, one can get a control on the quantities defined 
in \eqref{sup}. However, in many other cases of random iterates, such a control 
is not possible, while one can get some upper bounds on
\begin{equation}\label{L1}
\iint
\E \left (|X_{k, x}-X_{k,  y}|\right)  \nu(dx) \nu(dy)\, ,
\end{equation}
where $\nu$ is the invariant distribution of the chain 
$(W_n)_{n \geq 1}$. 

Consequently, we shall  establish two  distinct results, 
with different range of applicability. In Theorem  \ref{KMTavecdeltaninfty}, we give a strong approximation result under   conditions involving 
some quantities similar  to   \eqref{sup}. In Theorem 
\ref{KMTavecdeltan1}
the conditions are expressed in terms of the quantities
\eqref{L1}. 
The second Theorem applies to a large variety of examples,
including some  well known  examples of irreducible and
aperiodic Markov Chains with countable or continuous state space. 
These examples of ergodic Markov chains will allow us to prove  that the conditions 
given in Theorem \ref{KMTavecdeltan1}
are in some sense optimal. 

\smallskip

In all the paper, we shall use  the  notation
 $a_n \ll b_n$, which means that there exists a positive  constant $C$ not
depending on $n$ such that  $a_n \leq  Cb_n $, for all positive integers  $n$.

\section{Main results} \label{Sec2}

Let $(\Omega, {\mathcal A}, {\mathbb P})$ be a probability space, and let $(\varepsilon_i)_{i \geq 1}$ be iid random variables defined on $\Omega$, with values in a measurable space $G$ and with common distribution $\mu$. Let $W_0$ be a random variable defined on $\Omega$ with values in a measurable space $X$, independent of 
$(\varepsilon_i)_{i \geq 1}$, and let 
$F$ be a measurable function from $G\times X $ to $X$. For any $n \geq 1$, define
\[
W_{n} = F ( \varepsilon_{n}, W_{n-1} ) \, ,
\]
and assume that $(W_n, n \geq 1)$ has a stationary distribution $\nu$. 
Let now $h$ be a measurable function from $G\times X $ to ${\mathbb R}$ and define, for any $n \geq 1$,
\beq \label{functionofarandomiterates}
X_n = h ( \varepsilon_{n}, W_{n-1} ) \, .
\eeq
Then $(X_n)_{n \geq 1}$ forms a stationary sequence with stationary distribution, say  $\pi$.  Let $( {\mathcal G}_{i} )_{i \in {\mathbb Z}}$ be the non-decreasing filtration defined as follows: for any $i < 0$, 
${\mathcal G}_{i} =\{ \emptyset, 
\Omega\}$, ${\mathcal G}_{0} = \sigma (W_0)$ and for any $i \geq 1$, $ 
{\mathcal G}_{i} = \sigma ( \varepsilon_{i }, \ldots,  \varepsilon_{ 1 }, W_0 ) $. It follows that for any $n \geq 1$, $X_n$ is ${\mathcal G}_{n}$-measurable. 

%
%

\smallskip
Our first result proves that the strong approximation result holds with rate $n^{1/p}$ when the stationary distribution $\pi$  has a moment of order $p>2$ and we impose that the sequence of coupling coefficients $({\delta}_{\infty} (n) )_{n \geq 1}$ defined in \eqref{defdeltaninfty} decreases arithmetically to zero plus the condition \eqref{condonnorminftycarre}. As we shall see in Section \ref{sectionexamples}, these conditions are satisfied for instance  for the left random walk on $GL_d ( {\mathbb R})$. 

Let $W_0$ and $W_0^*$ be random variables with law $\nu$, and such that $W_0^*$ is independent of $(W_0, (\varepsilon_i)_{i \geq 1})$. For any $n \geq 1$,  let 
\beq \label{defXnstartWnstar}
X^*_n = h ( \varepsilon_{n}, W^*_{n-1} )\,  \mbox{ with } \,  W^*_{n} = F ( \varepsilon_{n}, W^*_{n-1} ) \, .
\eeq
Define then 
\beq \label{defdeltaninfty}
{\delta}_{\infty} (n)  = \Vert \E (|X_n - X_n^*| \, | (W_0,W_0^*)) \Vert_{\infty}  \, , n \geq 1 \, ,
\eeq 
where, above and in all the rest of the paper, the infinite norm is the usual essential supremum norm. 
\begin{Theorem} \label{KMTavecdeltaninfty} Let $(X_n, n \geq 1)$ be the stationary sequence defined by \eqref{functionofarandomiterates} and assume that its stationary distribution $\pi$  has moment of order $p>2$. Assume in addition that there exists a positive constant $c$ such that for any $n \geq 1$, 
\beq \label{condondeltaninfty}
{\delta}_{\infty} (n)
\leq  c n^{-q} \,  
\text{ with $q >(p-1)/2$} \, ,
\eeq 
where $(\delta_{\infty} (n))_{n \geq 1}$ is defined in \eqref{defdeltaninfty}, and that 
\beq \label{condonnorminftycarre}
\sup_{n \geq 1} \Vert \E (X_n^2 | {\mathcal G}_{n-1}) \Vert_{\infty} \leq  c \, .
\eeq 
Let $S_n = \sum_{k=1}^n X_k$. Then $n^{-1} \E \big ( (S_n - n \E (X_1) )^2 \big ) \rightarrow \sigma^2$ as $n \rightarrow \infty$ and one can redefine $(X_n)_{n \geq 1}$ without changing its distribution on a (richer) probability space on which 
there exist iid random variables $(N_i)_{i \geq 1}$ with common distribution ${\mathcal N} (0, \sigma^2)$, such that,
\beq \label{resultSA}
 S_n - n \E (X_1)- \sum_{i=1}^n N_i = o(n^{1/p}) \, \text{ ${\mathbb P}$-a.s.}
\eeq
\end{Theorem}

In the rest of this section, we shall give conditions expressed in terms of the quantities $\Vert X_n - X_n^* \Vert_1$ for the strong approximation \eqref{resultSA} to hold.  Before stating the result, we need to introduce some notations: 
\smallskip

For any $n \geq 0$, let  us define the sequence $(\delta (n) )_{n \geq 0}$ as follows
\[
\delta (0) = \delta (1) = \E (|X_1|) \,  \text{ and } \, \delta (n) = 2^{-1} \sup_{k \geq n-1} \Vert X_k - X_k^* \Vert_1 \, , \, n \geq 2 \, . 
\]
These quantities are finite  if $\pi$ has a moment of order $1$. 

\smallskip
For any $x \geq 0$, denote by 
\[
\delta (x) = \delta ([x]) 
\]
and, for any $u \in [0,  \E (|X_1|) ]$,  let
\[
\delta^{-1} (u) =  \inf \{q \in {\mathbb N} \, : \, \delta (q) \leq u \} = \sum_{n \geq 0} {\bf 1}_{u < \delta (n) }  \, .
\]
Denote also by $Q$ the quantile function associated with $|X|$ where $X$ is a random variable with law $\pi$: it is then the generalized inverse of the tail function 
$t \mapsto \p (|X|  >t )=\pi((-\infty, -t[)+ \pi(]t, \infty))$. Let $H $ be the function from $[0,1]$ 
to $ {\mathbb R}^+$ defined by $H(x) = \int_0^x Q(u) du$. We shall assume the following condition
\beq \label{condonQdelta-2}
\sum_{n \geq 1} n^{p-2} \int_{0}^{\delta (n)} Q^{p-1} \circ H^{-1} (u) du < \infty \, .
\eeq
\begin{Theorem} \label{KMTavecdeltan1} Let $(X_n, n \geq 1)$ be a stationary sequence defined by \eqref{functionofarandomiterates} and assume that its stationary distribution $\pi$  has a moment of order $p>2$. Assume in addition that condition \eqref{condonQdelta-2} holds. Let $S_n = \sum_{k=1}^n X_k$. Then $n^{-1} \E \big ( (S_n - n \E (X_1) )^2 \big ) \rightarrow \sigma^2$ as $n \rightarrow \infty$ and    one can redefine $(X_n)_{n \geq 1}$ without changing its distribution on a (richer) probability space on which 
there exist iid random variables $(N_i)_{i \geq 1}$ with common distribution ${\mathcal N} (0, \sigma^2)$, such that,
\[
 S_n - n \E (X_1)- \sum_{i=1}^n N_i = o(n^{1/p}) \, \text{ ${\mathbb P}$-a.s.}
\]
\end{Theorem}
\begin{remark} \label{remarkequivcond} If we define 
\beq \label{defgammaetinvgamma}
\gamma (x) = H^{-1} (\delta ([x])) \mbox{ for any $x \geq 0$  and } \gamma^{-1} (u) =\delta^{-1} \circ H(u)  \mbox{ for any $u \in [0,1]$}  \, ,
\eeq
then condition \eqref{condonQdelta-2} can be rewritten as 
\beq \label{condonQdelta}
\sum_{n \geq 1} n^{p-2} \int_{0}^{\gamma (k)} Q^p (u) du < \infty \, ,
\eeq
which also reads as  
\beq \label{notaforR}
\int_0^1 R^{p-1} (u) Q(u) du < \infty \, \text{ where } \, R(u) = \gamma^{-1} (u)  Q(u) \, , 
\eeq
\end{remark}
\begin{remark} \label{remarksufficient} Sufficient conditions for  \eqref{condonQdelta-2} to hold in terms of moments (or weak moments) of $\pi$ can be given by using Lemma 2 in Dedecker and Doukhan \cite{DD03}. For instance, if 
\beq \label{condondelta}
\Vert X_1 \Vert_r  \, \mbox{ for some $r >p$, and } \, \sum_{n\geq 1} n^{(pr-2r+1)/(r-p)}\delta(n) <\infty \, , 
\eeq
then condition \eqref{condonQdelta-2} is satisfied. Note that in the case where $\Vert X_1 \Vert_{\infty}< \infty$, condition \eqref{condonQdelta-2}  is equivalent to  $\sum_{n\geq 1} n^{ p-2}\delta(n) <\infty$. 
\end{remark}

If we define the following meeting time 
\beq \label{definitionmeetingtime} 
T^* = \inf \{ k \in {\mathbb N} \, : \,  W_k = W_k^* \} \,  ,
\eeq
it follows that, for any $n \geq 2$, 
\[
\delta (n) \leq  \int_{0}^{\p_{\nu \otimes \nu} (T^* \geq n)} Q (u) du \, .
\]
Therefore the following corollary holds. 
\begin{corollary} \label{corKMTavecdeltan1} Let $(X_n, n \geq 1)$ be the stationary sequence defined by \eqref{functionofarandomiterates} and assume that its stationary distribution $\pi$  has a moment of order $p>2$. Assume in addition that 
\beq \label{condavecTstar}
\sum_{ n \geq 0 } (n+1)^{p-2} \int_{0}^{\p_{\nu \otimes \nu} (T^* \geq  n )} Q^p (u) du < \infty \, .
\eeq
Then the conclusions of Theorem \ref{KMTavecdeltan1} hold. 
 \end{corollary}
According to the computations given in Annex C of Rio \cite{Ri00}, if 
\beq \label{condontailofTstar}
\Vert X_1 \Vert_r  \, \mbox{ for some $r >p$, and } \, \sum_{n\geq 1} n^{(pr-2r+p)/(r-p)} \p_{\nu \otimes \nu} (T^* \geq  n )  <\infty \, , 
\eeq
then condition \eqref{condavecTstar} is satisfied. In the case where $\Vert X_1 \Vert_{\infty}< \infty$, condition \eqref{condavecTstar} is equivalent to 
\beq \label{condavecTstarbounded}
\sum_{ n \geq 1 } n^{p-2} \p_{\nu \otimes \nu} (T^* \geq  n )  < \infty \, .
\eeq
Propositions  \ref{propCE} and \ref{propCEcontinuous} in Section \ref{sectionEMC} will show that condition \eqref{condavecTstarbounded} is optimal in some sense.

\section{Applications} \label{sectionexamples}

\subsection{Left random 
walk on $GL_d({\mathbb R})$}\label{LRW}
As in the introduction, let $(\varepsilon_n)_{n \geq 1}$ be independent random matrices taking values in $G= GL_d(\mathbb R)$, $d \geq 2$, with common distribution $\mu$. 
let $A_0={\rm Id}$ and for every $n \geq 1$, $A_n = \varepsilon_n \cdots \varepsilon_1$.

Let $\Vert \cdot \Vert$ be the euclidean norm on ${\mathbb R}^d$.
Recall that  $\mu $ has a moment of order $p \geq 1$ if 
\eqref{Mpmu} holds. 
 Recall also that if $\mu$ admits a moment of order $1$ then 
   \eqref{LLN} holds, and the quantity $\lambda_\mu$ is well defined.

Let $X:=P_{d-1} ({\mathbb R}^d)$ be the projective space of ${\mathbb R}^d -\{0\}$ and write ${\bar x}$ as the projection of $x \in {\mathbb R}^d -\{0\}$ to $X$. We assume that $\mu$ is strongly irreducible (i.e. that no proper finite union of subspaces of $\R^d$ are invariant by $\Gamma_\mu$, the closed semi-group generated by the support of $\mu$) and proximal (i.e. that there exists a matrix in $\Gamma_\mu$ 
admitting a unique (with multiplicity one) eigenvalue with maximum modulus). 
Under those assumptions (see e.g. Bougerol-Lacroix \cite{BL} or 
Benoist-Quint \cite{BQ}) it is well-known that there exists a unique invariant measure $\nu$ on ${\mathcal B} (X)$, meaning that for any continuous and bounded function $f$ from $X$ to $\mathbb R$,
\[
\int_X f(x) \nu(dx) = \int_G \int_X f( g \cdot x ) \mu(dg) \nu(dx) \, .
\]
The left random walk of law $\mu$ is the process defined by $W_0:=\varepsilon_0$ and $W_n =\varepsilon_n W_{n-1}$ for $n \geq 1$ where we assume that $\varepsilon_0$ is independent of $(\varepsilon_n)_{n \geq 1}$. 
As explained in the introduction, our aim is to study the partial sums associated with the 
random sequence   $(X_n)_{n \geq 1}$ given by 
\[
X_n := h (\varepsilon_n, W_{n-1} ) \, , \, n \geq 1 \, , 
\]
where for every $g \in G$ and every ${\bar x} \in X$,
\[
 h ( g , {\bar x} ) = \log \Big ( \frac{\Vert g \cdot x \Vert }{ \Vert x \Vert }\Big ) \, .
\]

As usual, we shall denote by $X_{n,{\bar x}}$ the random variable for which $W_0 = {\bar x}$. We then define $S_{n,{\bar x}}= \sum_{k=1}^n X_{n,{\bar x}}$ and recall that the identity \eqref{mainidentity} holds: for any $x \in S^{d-1}$, 
\begin{equation*}
S_{n, \bar x} = \sum_{k=1}^n X_{k, \bar x} = \log \Vert  A_n x \Vert \, .
\end{equation*}
Applying Theorem  \ref{KMTavecdeltaninfty}, the following strong  approximation with rate holds. 
\begin{corollary} \label{KMTforcocycles} Let $\mu$ be a proximal and strongly irreducible probability measure on ${\mathcal B}(G)$. Assume that $\mu$ has a moment of order $p>2$. Then $n^{-1} \E_{\nu} \big ( (S_n - n \lambda_{\mu} )^2\big ) \rightarrow \sigma^2$ as $n \rightarrow \infty$ and for every (fixed) ${\bar x} \in X$, one can redefine $(S_{n,{\bar x}})_{n \geq 1}$ without changing its distribution on a (richer) probability space on which 
there exist iid random variables $(N_i)_{i \geq 1}$ with common distribution ${\mathcal N} (0, \sigma^2)$, such that,
\[
 S_{n,{\bar x}} - n \lambda_{\mu} - \sum_{i=1}^n N_i = o(n^{1/p}) \, \text{a.s.}
\]
\end{corollary}
\begin{remark} It follows from item $c)$ of Theorem 4.11 of 
Benoist-Quint \cite{BQ} that $\sigma>0$ if $\mu$ is strongly irreducible and the image of $\Gamma_\mu$ 
in $PGL_d(\R)$ is unbounded.  
\end{remark}

\noindent {\bf Proof of Corollary \ref{KMTforcocycles}.} Using the same arguments as in Cuny-Dedecker-Jan \cite{CDJ}
(see the proof of their Theorem 1), we infer that it suffices to prove the result on stationary regime. More precisely, it suffices to prove that  one can redefine $(S_{n})_{n \geq 1}$ without changing its distribution on a (richer) probability space on which 
there exist iid random variables $(N_i)_{i \geq 1}$ with common distribution ${\mathcal N} (0, \sigma^2)$, such that,
\beq \label{KMTSR}
 S_{n} - n \lambda_{\mu} - \sum_{i=1}^n N_i = o(n^{1/p}) \, \text{ ${\mathbb P}_{\nu}$-a.s.}
\eeq
Note also that the fact that $n^{-1} \E_{\nu} \big ( (S_n - n \lambda_{\mu} )^2\big ) \rightarrow \sigma^2$ as $n \rightarrow \infty$ comes from Theorem 2 (ii) in  \cite{CDJ}. Now the strong invariance principle  \eqref{KMTSR} is a direct application of Theorem  \ref{KMTavecdeltaninfty}. To see this, note first that the following estimate is valid (see Proposition 3 in  \cite{CDJ}):
\[
 \sum_{k\geq 1}  k^{p-2} \sup_{{\bar x }, {\bar y } \in X} \E \big (  \big | X_{k,{{\bar x}}} - X_{k,{{\bar y}}} \big | \big ) < \infty  \, .
\]
Since $\big(\sup_{{\bar x }, {\bar y } \in X} \E \big (  \big | X_{k,{{\bar x}}} - X_{k,{{\bar y}}} \big |\big)_{k\ge 1}$ is non
 increasing,   $\sup_{{\bar x }, {\bar y } \in X} \E \big (  \big | X_{k,{{\bar x}}} - X_{k,{{\bar y}}} \big |  \big ) \ll k^{- (p-1)}$. Hence condition \eqref{condondeltaninfty}  holds with $q =p-1$. To end the proof it suffices to notice that 
condition \eqref{condonnorminftycarre} also holds since, for any $k \geq 1$, 
\[\left \| {\mathbb E}  ( X_{k}^2  | 
\mathcal{G}_{k-1} ) \right \|_\infty \leq 
\int_G  (\log N(g))^2 \mu(dg) < \infty \, .
\]
\hfill $\square$

\subsection{Contracting iterated random functions}

\subsubsection{Uniform contraction}

Assume that there is a distance $d$ on $X$, and that there
exist $\kappa>0$ and $\rho \in (0,1)$ such that, for any $n \geq 1$, 
\begin{equation}\label{contr}
\Vert \E (d(W_n , W_n^*) \, | (W_0,W_0^*)) \Vert_{\infty} \leq \kappa \rho^n  \, ,
\end{equation}
where $W_n^*$ is defined in \eqref{defXnstartWnstar}. Note that condition \eqref{contr} holds  if the chain is 
``one step contracting" in the following sense
\[
\|d(W_0,W_0^*)\|_{\infty} < \infty \  \text{and} \ 
{\mathbb E} (d(W_{1,x}, W_{1,y}))\leq \rho d(x,y) 
\ \text{for any $(x,y) \in X \times X$} \, .
\]

Let us now define a class of observables from $G\times X$ to 
${\mathbb R}$ 
 for which 
one can easily compute the coefficient $\delta_\infty(n)$. 
Let 
$\eta$ be a measurable function from $G$ to ${\mathbb R}^+$
such that ${\mathbb E}(\eta(\varepsilon_0)) < \infty$, and 
let $c$ be a concave non-decreasing function from ${\mathbb 
R}^+$ to ${\mathbb R}^+$ such that $c(0)=0$. 

One says that $h :  G \times X \rightarrow \R$ belongs to the 
class ${\mathcal L}(\eta, c)$ if, 
$$
 |h(z, x)-h(z,y)|
 \leq  \eta(z) c(d(x,y)) \ \text{for any $(x,y,z) \in X \times X \times G$}\, .
$$

\begin{Lemma}\label{fastoche}
Assume that the stationary Markov chain $(W_n)_{n \geq 0}$ satisfies the contraction condition \eqref{contr}, and
let  $(X_n)_{n \geq 1}$ be defined by 
\eqref{functionofarandomiterates} for some $h \in {\mathcal L}(\eta, c)$. Then, there exists a constant $A>0$ such that, for any $n \geq 1$, 
$$
\delta_{\infty}(n) \leq A  c(\kappa \rho^{n-1}) \, .
$$
\end{Lemma}

\noindent {\bf Proof.} 
Let $A={\mathbb E}(\eta(\varepsilon_0))$.  Since $h$ belongs to ${\mathcal L}(\eta, c)$, and since 
$c$ is concave, 
$$
\E (|X_{n+1,x} - X_{n+1,y}| ) \leq A \E(c(d(W_{n,x}, W_{n,y}))) \leq A c( \E(d(W_{n,x}, W_{n,y})))\, .
$$
Hence, since $c$ is non-decreasing and $(W_n)_{n \geq 0}$
satisfies \eqref{contr},
\[
\Vert \E (| X_{n+1} - X_{n+1}^*| \, | (W_0,W_0^*)) \Vert_{\infty} 
  \leq A c( \Vert \E (d(W_n , W_n^*) \, | (W_0,W_0^*)) \Vert_{\infty} )
\leq A c( \kappa \rho^n)\, .   
\]
\hfill $\square$

\smallskip

Applying Theorem \ref{KMTavecdeltaninfty}, the following result holds:

\begin{corollary}\label{corcontractinfini}
Assume that the stationary Markov chain $(W_n)_{n \geq 0}$ satisfies the contraction condition \eqref{contr}, and
let  $(X_n)_{n \geq 1}$ be defined by 
\eqref{functionofarandomiterates} for some $h \in {\mathcal L}(\eta, c)$. Assume moreover that ${\mathbb E}(\eta( \varepsilon_1)^p) < \infty$ for some $p>2$, and that there
exists
 $x_0 \in X$ such 
that $\|c(d(W_0, x_0))\|_\infty < \infty$ and 
${\mathbb E}(|h(\varepsilon_1, x_0)|^p) < \infty$. 
If 
$c(\kappa \rho^n)=  O( n^{-q})$  
for some  $q >(p-1)/2$, then the conclusion of Theorem 
\ref{KMTavecdeltaninfty}
holds.
\end{corollary}

\begin{remark}
Note that Corollary \ref{corcontractinfini} applies to 
a large class of continuous observales (as functions of $x$), including all H\"older observables 
(case where $c(x)=x^\alpha$ for some $\alpha \in (0,1)$). More precisely it applies to any concave non-decreasing function  $c$ such 
that $c(x) \leq C |\ln (x)|^{-\gamma}$ in a neighborhood
of $0$, 
for some $\gamma > (p-1)/2$. 
\end{remark}

\noindent {\bf Proof of Corollary \ref{corcontractinfini}.}  Applying  Lemma \ref{fastoche}, we infer that $\delta_\infty(n)=  O( n^{-q})$  
for some  $q >(p-1)/2$. Hence, if one can prove that 
\begin{equation}\label{Mp}
\sup_{n \geq 1} \Vert \E (|X_n|^p | {\mathcal G}_{n-1}) \Vert_{\infty} \leq  M \, , 
\end{equation}
for some finite constant $M$, the result will follow directly from 
Theorem \ref{KMTavecdeltaninfty}. 
 To prove 
\eqref{Mp}, we note that 
\begin{equation}\label{dec}
\E (|h(\varepsilon_n, W_{n-1})|^p | {\mathcal G}_{n-1})
\leq 2^{p-1} 
\E (|h(\varepsilon_n, W_{n-1})- h(\varepsilon_n, x_0)|^p | {\mathcal G}_{n-1})
+
2^{p-1} {\mathbb E}(|h(\varepsilon_1, x_0)|^p) \, .
\end{equation}
For the first term on the right-hand side of \eqref{dec}, we use the fact that $h \in {\mathcal L}(\eta, c)$, which gives 
\begin{equation}\label{last}
\E (|h(\varepsilon_n, W_{n-1})- h(\varepsilon_n, x_0)|^p | {\mathcal G}_{n-1}) \leq 
{\mathbb E}(\eta( \varepsilon_1)^p) \|c(d(W_0, x_0))\|_\infty^p \, .
\end{equation}
Under the assumptions of Corollary \ref{corcontractinfini}, it follows from 
\eqref{dec} and \eqref{last} that the upper bound 
\eqref{Mp} holds. \hfill $\square$

\subsubsection{${\mathbb L}^1$-contraction}

Assume that there is a distance $d$ on $X$, and that there
exist $\kappa>0$ and $\rho \in (0,1)$ such that, for any $n \geq 1$, 
\begin{equation}\label{contr1}
\|d(W_{n}, W_{n}^*)\|_1 \leq \kappa \rho^n  \, ,
\end{equation}
where $W_n^*$ is defined in \eqref{defXnstartWnstar}.
Note that condition \eqref{contr1} holds  if the chain is 
``one step contracting" in the following sense: 
\[{\mathbb E} (d(x_0, F(\varepsilon_1, x_0)) < \infty   \  \text{for some $x_0 \in X$}
\]
and
\[
{\mathbb E} (d(W_{1,x}, W_{1,y}))\leq \rho d(x,y) 
\ \text{for any $(x,y) \in X \times X$} \, .
\]
Note also that, under the two conditions above, there exists an unique  stationary distribution $\nu$ (see Theorem 2 of \cite{SW04}).

Let us now define a class of observables from $G\times X$ to 
${\mathbb R}$ 
 for which 
one can easily compute the coefficients $\delta (n) $. 
Let $c$ be a concave non-decreasing function from ${\mathbb 
R}^+$ to ${\mathbb R}^+$ such that $c(0)=0$. 

One says that $h :  G \times X \rightarrow \R$ belongs to the 
class ${\mathcal L}(c)$ if, 
\[
\E(|h(\varepsilon_1, x)-h(\varepsilon_1,y)|)
 \leq   c(d(x,y)) \ \text{for any $(x,y) \in X \times X $}\, .
\]

\begin{Lemma}\label{fastoche1}
Assume that the stationary Markov chain $(W_n)_{n \geq 0}$ satisfies the contraction condition \eqref{contr1}, and
let  $(X_n)_{n \geq 1}$ be defined by 
\eqref{functionofarandomiterates} for some $h \in {\mathcal L}(c)$. Then, for $n \geq 2$, 
\[
\delta(n) \leq  2^{-1} c(\kappa \rho^{n-2}) \, .
\]
\end{Lemma}

\noindent {\bf Proof.} 
Let $k \geq n \geq 2$.  Since $h$ belongs to ${\mathcal L}(c)$, and since 
$c$ is concave, 
$$
\|X_{k} - X^*_{k}\|_1\leq  \|c(d(W_{k-1}, W^*_{k-1}))\|_1 \leq  c\left (  \|d(W_{k-1}, W^*_{k-1})\|_1 \right ) \, .
$$
Hence, since $c$ is non-decreasing and $(W_n)_{n \geq 0}$
satisfies \eqref{contr1},
\[
 \|X_{k} - X^*_{k}\|_1  \leq 
 c\left ( \kappa \rho^{k-1} \right )\, .   
\]
The result follows from the definition of $\delta(n)$
and the fact that $c$ is non-decreasing. \hfill $\square$

\smallskip

Recall that the function $Q$ and $H$ related to 
the tail function $t \mapsto {\mathbb P}(|X_1|>t)$
have been defined in Section \ref{Sec2}.
Combining Theorem \ref{KMTavecdeltan1} and Lemma \ref{fastoche1}, the following result  holds:

\begin{corollary}\label{corcontract1}
Assume that the stationary Markov chain $(W_n)_{n \geq 0}$ satisfies the contraction condition \eqref{contr1}, and
let  $(X_n)_{n \geq 1}$ be defined by 
\eqref{functionofarandomiterates} for some $h \in {\mathcal L}(c)$. Assume moreover that
\begin{equation}\label{cond1c}
\sum_{n \geq 1} n^{p-2} \int_{0}^{c\left (\kappa \rho^n \right )} Q^{p-1} \circ H^{-1} (u) du < \infty \, .
\end{equation}  Then the conclusion of Theorem 
\ref{KMTavecdeltan1}
holds.
\end{corollary}

\begin{remark}
From Remark \ref{remarksufficient}, it follows that \eqref{cond1c} holds as 
soon as
\begin{equation} \label{condonc}
\Vert X_1 \Vert_r  \, \mbox{ for some $r >p$, and } \, \sum_{n\geq 1} n^{(pr-2r+1)/(r-p)}c\left (\kappa \rho^n \right ) <\infty \, . 
\end{equation}
The condition \eqref{condonc} is equivalent to the following integral condition on the function $c$
$$
  \int_0^{1/2} \frac{1}{t}c(t) |\ln(t)|^{(pr-2r+1)/(r-p)}\, dt < \infty \, .
$$
\end{remark}

\subsection{Ergodic Markov chains} \label{sectionEMC}

\subsubsection{A discrete ergodic Markov chain example}
Let  $(\varepsilon_{i})_{i \in \Z}$ be a sequence of iid real-valued random variables distributed as $\varepsilon $ with 
\[
\p (\varepsilon = k ) =p_k \, , \, k \in {\mathbb N^*} \, ,
\]
Let $W_0$ be a random variable with values in $\N$ independent of $(\varepsilon_{i})_{i \in \Z}$, and define for any $k \geq 1$,
\beq \label{defWkDMC}
W_k = (W_{k-1} -1) {\bf 1}_{W_{k-1} \neq 0} + (\varepsilon _k -1) {\bf 1}_{W_{k-1} =0 } \, .
\eeq
Hence $(W_k , k \in {\mathbb N}) $ is a Markov chain with state space ${\mathbb N}$, initial distribution ${\mathcal L} (W_0)$ and transition probabilities satisfying
\[
P_{i,i-1} =1  \mbox{ and } P_{0,i-1} =p_i \, , \, i \geq 1 \, .
\]
Assume that $p_1>0$ and $p_{n_j} >0$ along $n_j  \rightarrow \infty$. Then the chain  $\{W_{k};k\geq 0\}$ is irreducible and aperiodic. Moreover, the stationary distribution exists 
if and only if $\E(\varepsilon )<\infty $ and  is given by
\begin{equation*}
\nu _{0}=1/\E(\varepsilon ) \,  \mbox{ and } \, \nu _{j}=\nu _{0}\sum_{i=j+1}^{\infty }p_{i}\;,\;j \in {\mathbb N}^*\, .
\end{equation*}
\begin{corollary} \label{cordiscreteMC} Let $ p >2$ and $f$ be a function from ${\mathbb N} $ to ${\mathbb R}$ such that $\nu (|f|^r) < \infty$ with $r >p$. Assume that
\beq \label{condMCdiscrete}
\sum_{n \geq 1} n^{\frac{p (r-1)}{r-p}} p_n < \infty \, .
\eeq
Then condition \eqref{condontailofTstar} is satisfied and the conclusions of Theorem \ref{KMTavecdeltan1} hold  for  $X_n = f ( W_n)$ where $(W_n)_{n \geq 0}$ is the Markov chain defined by \eqref{defWkDMC} with ${\mathcal L} (W_0)=\nu$.
\end{corollary}
For bounded observables (case $r = \infty$), condition  \eqref{condMCdiscrete} reads as $\sum_{n \geq 1} n^{p} p_n < \infty$.  As we shall see in the proof of the next proposition (see \eqref{toshowtheequiavalence}),  $\sum_{n \geq 1} n^{p} p_n < \infty$ is equivalent to $\sum_{ n \geq 1 } n^{p-2} \p_{\nu \otimes \nu} (T^* \geq n ) < \infty$, where  $T^*$ is the meeting time defined in \eqref{definitionmeetingtime}. The next proposition shows that this latter condition is in some sense optimal. 
\begin{Proposition} \label{propCE}
Let $p >2$ and $(W_k)_{k \geq 0}$ be the Markov chain described above with $p_k: = 1/(\zeta(p+1)k^{p+1})$, $k \in {\mathbb N^*} $,  where 
$\zeta(p+1) = \sum_{k \geq 1} k^{-(p+1)} $. Then 
\beq \label{cond1propCE}
\sum_{ n \geq 1 } n^{p-2} \p_{\nu \otimes \nu} (T^* \geq n )  = \infty  \text{, and for any $\varepsilon >0$, } \, \sum_{ n \geq 2 } n^{p-2} (\log n)^{-(1+ \varepsilon)} \p_{\nu \otimes \nu} (T^* \geq n )  < \infty \, .
\eeq
Moreover, for any stationary and Gaussian centered sequence $(g_k)_{k \in {\mathbb Z}}$ with  convergent
series of covariances,
\beq \label{cond2propCE0}
\limsup_{n \rightarrow \infty}  \, (n \log n)^{-1/p} \Big |  \sum_{k=1}^n ( {\mathbf 1}_{\{W_k =0\}} - \nu_0) -\sum_{k=1}^n g_k \Big | >0 \mbox{ almost surely.}
\eeq
\end{Proposition}
\noindent {\bf Proof of Corollary \ref{cordiscreteMC}.}  Define
\[
T^*_0 = \inf \{ k \in {\mathbb N} \, : \,  W_{k} = W^*_{k} =0  \}  \, .
\]
By definition,  $T^* \leq T_0^*$.  Hence for any $n \in {\mathbb N}$,
\beq \label{relation1TTstar}
\p_{\nu \otimes \nu}  ( T^* \geq n ) \leq \p_{\nu \otimes \nu}  ( T_0^* \geq n ) \, .
\eeq
Next, it is easy to see that  for any $n \in {\mathbb N}$,
\beq \label{equalityoflaw}
\p_{\nu \otimes \nu}  ( T_0^* = n  ) = \p_{\nu \otimes \nu}  ( T'_0= n  )
\eeq
with 
\[
T'_0 = \inf \{ k \geq 0 \, : \,  W_{k} = W_k' =0  \}  \, ,
\]
where $(W_k' , k \in {\mathbb N}) $ is the Markov chain defined as follows: Let  $(\varepsilon'_k)_{k \in {\mathbb Z}}$ be an independent copy of  $(\varepsilon_k)_{k \in {\mathbb Z}}$ and 
independent of $W_0$. Let $W'_0$ be independent of $(W_0, (\varepsilon_k)_{k \in {\mathbb Z}}, (\varepsilon'_k)_{k \in {\mathbb Z}})$ and, 
for any $k \geq 1$, set
\[
W'_k = (W'_{k-1} -1) {\bf 1}_{W'_{k-1} \neq 0} + \varepsilon' _k  {\bf 1}_{W'_{k-1} =0 } \, .
\]
According to Lindvall \cite{Li79}, if $\E_{\nu} (\psi (\tau)) < \infty$ where $\tau = \inf \{ k \geq 1 \, : \, W_{k} =0 \}$ and $\psi $ is a non-decreasing function from ${\mathbb N}$ to $[2, \infty [$ such that 
$( ( \log (\psi (n))/n)_n$ is non-increasing and converges to $0$, then $\E_{\nu \otimes \nu} (\psi (T'_0)) < \infty$. Note now that 
\beq \label{boundtailtauCMdiscrete}
\p_{\nu} ( \tau \geq n ) = \sum_{\ell =0}^{n-1} \nu_{\ell} \sum_{j \geq n- \ell} p_j + \sum_{\ell \geq n } \nu_{\ell } \leq  2 \sum_{\ell \geq [n/2] } \nu_{\ell } + \sum_{j \geq [n/2] +1} p_j\, .
\eeq
Hence under \eqref{condMCdiscrete},  $\E_{\nu} (\psi_{r,p} (\tau)) < \infty$ with $\psi_{r,p} (x) = x^{ r (p  -1) / (r-p)}$.  It follows that   $\E_{\nu} (\psi_{r,p} (T'_0)) < \infty$ which in turn implies that $\E_{\nu} (\psi_{r,p} (T^*)) < \infty$ by taking into account \eqref{relation1TTstar} and \eqref{equalityoflaw}. Therefore condition \eqref{condontailofTstar} is satisfied and Corollary \ref{corKMTavecdeltan1} applies. \hfill $\square$

\medskip

\noindent {\bf Proof of Proposition \ref{propCE}.} Note first that the following coupling inequality holds: for any $n \geq 1$, 
\beq \label{couplinginebeta}
\beta (n) :=  \frac{1}{2}  \int \Vert  \delta_x P^n  - \nu  \Vert_{v} \nu (dx)  \leq  \p_{\nu \otimes \nu}  ( T^* \geq  n ) \, , 
\eeq
where $\Vert \mu \Vert_{v}$ denotes the total variation norm of a signed measure $\mu$ and $P$ is the transition function of the Markov chain $(W_k)_ {k \in {\mathbb N}}$.  But for any $n \geq 1$, $\beta (n) \geq 2 \alpha (n)$ where $(\alpha (n))_{n \geq 1}$ is the sequence of strong mixing coefficients of the chain which starts from the stationary distribution. As quoted in Chapter 30 of Bradley \cite{Br07}, 
\[
 \alpha (n) \geq \big | \p_{\nu}(  W_0  \geq n+1 , W_n =0) - \p_{\nu}(  W_0  \geq n+1) \p_{\nu} (W_n =0) \big | = \nu_0 \sum_{k \geq n+1} \nu_k  \, .
\]
It follows that for any $s \geq 0$, 
\[
\sum_{n \geq 1} n^s \sum_{k \geq n+1} \nu_k = \infty \Rightarrow \sum_{n \geq 1} n^s \p_{\nu \otimes \nu}  ( T^* \geq n ) = \infty \, ,
\]
which together with the arguments developed in the proof of  Corollary \ref{cordiscreteMC} show that 
\beq \label{toshowtheequiavalence}
\sum_{n \geq 1} n^{s+2 } p_n  < \infty \iff \sum_{n \geq 1} n^s \p_{\nu \otimes \nu}  ( T^* \geq n ) < \infty \, .
\eeq
This proves the first part of \eqref{cond1propCE}. To prove its second part, it suffices to use again the arguments developed in the proof of  Corollary \ref{cordiscreteMC} and to notice 
that, for $p_k: = 1/(\zeta(p+1)k^{p+1})$, $k \in {\mathbb N^*} $, the upper bound \eqref{boundtailtauCMdiscrete} entails that  $\E_{\nu} (\psi_p (\tau)) < \infty$ with $\psi_p (x) = \frac{ x^{p-1}}{(\log (1+x))^{1 + \varepsilon }}$ where $\varepsilon >0$. This ends the proof of \eqref{cond1propCE}. 
\medskip

To prove the second part of the proposition, we shall use similar arguments as those developed in the proof of Theorem 2.2 in Dedecker-Merlev\`ede-Rio \cite{DMR14} and 
adopt the following notations:  the regeneration times $(T_k)_{k\geq 0}$ of the Markov chain $(W_k)_{ k \geq 0}$ are defined by induction as follows: 
$R_0 = \inf \{n>0 \, : \,  W_n =0 \}$ and $R_k = \inf \{ n > R_{k-1}  \, : \,  W_n =0 \}$. Let $\tau_k = R_{k+1} - R_k$ for $k \geq 0$. Note that $(\tau_k)_{k \geq 0}$ are iid and that their common law is the law of $R_0$ 
when the chain starts at zero. Note that 
\[
\p (\tau_k >  (k\ln k)^{1/p} ) \sim  c_p /(k\log k) \, .
\]
Since the regeneration times $\tau_k$ are independent, by the converse Borel-Cantelli lemma, it follows that 
\[
\p (\tau_k >  (k\log k)^{1/p}  \text{ infinitely often }) = 1 \, . 
\]
Now we take 
\[
f(x) := {\mathbf 1}_{x =0} \mbox{ and } \, g(x) := f(x) - \nu_0 \, .
\]
$f$ is obviously a bounded function and $\nu (g) = 0$. Note that, for any $\ell \geq 0$, 
\[
\sum_{k= R_{\ell} +1}^{ R_{\ell } + m } ( \nu_0/2 - f(W_k)) = m \nu_0/2 \mbox{ for any $1 \leq m < \tau_{\ell}$} \, .
\]
Since $R_n/n$ converges to $\bkE (\tau_0)$ almost surely, it follows that, for some positive constant $c$ depending on $\bkE (\tau_0)$, 
\begin{equation}  \label{incresums}
\limsup_n   \sum_{i=n+1}^{n+ [c (n\log n )^{1/p} ] } ( \nu_0/2 - f(W_i)  )\geq  0 \text{ almost surely}. 
\end{equation}
Consider now a stationary and Gaussian centered sequence $(g_k)_{k\in {\mathbb Z}}$ with convergent series of covariances.  If follows from  both the Borel-Cantelli lemma
and the usual tail inequality for Gaussian random variables that, for any positive $\theta$,  
\[
\liminf_n \sum_{i=n+1}^{n+ [c (n\log n )^{1/p} ] } (g_i +  \theta ) \geq  0 \text{ almost surely.}  
\]
Taking $\theta =  \nu_0/4$ in the above inequality and using \eqref{incresums}, we then infer  that 
$$
\limsup_{n\rightarrow \infty} \frac{1}{[c (n\log n )^{1/p} ] }  \sum_{i=n+1}^{n+ [c (n\log n )^{1/p} ] }  \Big (   g_i + \nu_0  -f(W_i) \Big ) \geq  \nu_0/4 \ \text{ almost surely,} 
$$
which implies \eqref{cond2propCE0}. \hfill $\square$

\subsubsection{An example of ergodic Markov chain with continuous state space}

In this section, we consider an homogenous Markov chain with state space $[0,1]$ and transition probability
kernel $P(x, \cdot) $ given by 
\beq  \label{transitionproba}
P(x,A)=(1-x)\delta_{x}(A)+x\pi(A)\, ,
\eeq
where $\delta_{x}$ denotes the Dirac measure  at point $x$ and
\[
\pi (dx) = (a+1) x^a  dx  \, \text { with  $a >1$.}
\]
Note that the chain is irreducible and aperiodic and admits a unique invariant probability measure $\nu$ given by 
\[
 \nu (dx) = a x^{a-1} dx  \, .
\]
As in Section 9.3 in Rio \cite{Ri00}, we now construct a stationary Markov chain $(W_n)_{n \in {\mathbb N}}$ with initial law $\nu$ and transition probability
measure $P(x, \cdot)$. Let $\xi_0$ be a random variable with law $\nu$. We assume that the underlying
probability space is rich enough to contain a sequence $(\varepsilon_i)_{i \in {\mathbb Z}}:= (U_i, V_i)_{i \in {\mathbb Z}}$ of independent random
variables with uniform law over $[0,1]^2$, and that this random sequence is independent of $\xi_0$. The stationary Markov chain $(W_n)_{n \in {\mathbb N}}$ is then constructed 
via the following recursive equation: $W_0 = \xi_0$ and, for any $k \geq 1$,
\beq  \label{ecritureMC}
W_k =  W_{k-1} {\bf 1}_{U_{k} \geq W_{k-1}} + F_{\pi}^{-1} (V _k)  {\bf 1}_{U_{k}  < W_{k-1} } \, ,
\eeq
where $F_{\pi}^{-1}$ is the inverse of the cumulative function of $\pi$. It is easy to see that $(W_n)_{n \in {\mathbb N}}$ is a Markov chain with initial distribution $\nu$ and transition probability
kernel given by \eqref{transitionproba}. 
\begin{corollary} \label{corcontinuousMC} Let $ p >2$ and  $(W_k)_{k \in {\mathbb N}}$ be the stationary Markov chain defined by \eqref{ecritureMC} with $a >p-1$. Then condition 
\eqref{condavecTstarbounded} is satisfied and the conclusions of Theorem \ref{KMTavecdeltan1} hold  for  $X_n = f ( W_n)$, for any bounded function $f$ defined on $[0,1]$. 
\end{corollary}
The proof of this corollary is a direct application of Corollary \ref{corKMTavecdeltan1} by taking into account  the  following lemma whose proof is postponed to the Appendix (see Section \ref{sectionprooflmacontinuous}). 
\begin{Lemma} \label{lmacontinuous} For any $a>1$ there exist  positive constants $c(a) $ and  $C(a)$ depending only on $a$ such that for any $n \geq 1$,
\beq \label{aimMC}
\frac{c(a)}{n^a}  \leq \p_{\nu \otimes \nu} (T^* > n) \leq \frac{C(a)}{n^a} \, ,
\eeq
where $T^*$ is the meeting time defined in \eqref{definitionmeetingtime}. 
\end{Lemma}
In addition, this lemma together with Theorem 2.2 in Dedecker-Merlev\`ede-Rio \cite{DMR14} proves the sharpness of condition 
\eqref{condavecTstarbounded}  also in case of Markov chains with continuous state space. This is summarized in the next proposition. 
\begin{Proposition} \label{propCEcontinuous}
Let $p >2$ and $(W_k)_{ k \in {\mathbb N}}$ be the stationary Markov chain defined by \eqref{ecritureMC} with $a =p-1$. Then condition 
\eqref{condavecTstarbounded} fails. In addition,  for any map $f$ from $[0,1]$ to ${\mathbb R}$ with continuous and strictly positive derivative $f'$ on $[0,1]$, and  any stationary and Gaussian centered sequence $(g_k)_{k \in {\mathbb Z}}$ with  convergent
series of covariances,
\beq \label{cond2propCE}
\limsup_{n \rightarrow \infty}  \, (n \log n)^{-1/p} \Big |  \sum_{k=1}^n  f(W_k^a ) - n \int_0^1 f(t) dt  -\sum_{k=1}^n g_k \Big | >0 \mbox{ almost surely.}
\eeq
\end{Proposition}

\subsection{Lipschitz autoregressive models} 

We consider the autoregressive Lipschitz model as in Dedecker-Rio \cite{DR00}.
Let $\tau\in[0,1)$, $C\in (0,1]$ and $f\, :\, \R\to \R$ a $1$-Lipschitz 
function such that 

$$
f(0)=0\quad \mbox{ and } \quad |f'(t)|\le 1-\frac{C}{(1+|t|)^\tau} 
\quad \mbox{ for almost every $t$} \, .
$$
Let $(\varepsilon_i)_{i\ge 1}$ be iid real-valued random valued with common law $\mu$ and define for any $n \geq 1$
\beq \label{defLAM}
W_n = f (W_{n-1} ) + \varepsilon_n \, , \, \mbox{ with $W_0$ independent of $(\varepsilon_i)_{i\ge 1}$.}
\eeq
Let $S_n (h) = \sum_{k=1}^n h(W_i)$ for any measurable function $h$. 

\smallskip

The model above corresponds to the previously considered situation with 
$G=\R$ and $F\, :\, \R\times \R \to \R$ given by 
$F(x,y)=x+f(y)$, for every $x,y\in \R$.

Let $S\ge 1$ and assume that $\mu$ admits a moment of order $S$. 
It follows from Dedecker-Rio \cite{DR00} that there exists a unique invariant 
probability $\nu$ on $\R$, such that 
\begin{equation}\label{DR-moment}
\int_\R|x|^{S-\tau}\nu(dx) <\infty\, .
\end{equation}

The following strong approximation with rates holds.

\begin{corollary}
Let $\tau\in (0,1)$ and assume that $\mu$ admits a moment of order $S=p+\tau p$ for some $p>2$. Let $(W_n)_{n \geq 0}$ be defined by \eqref{defLAM} with ${\mathcal L} (W_0) = \nu$. 
Then, for any Lipschitz function $h$ such that $\nu (h) =0$,  $n^{-1}{\rm Var} (S_n(h) ) \rightarrow \sigma^2 (h)$ as $n \rightarrow \infty$ and   
one can redefine $(W_n)_{n \geq 0}$ without changing its distribution on a (richer) probability space on which 
there exist iid random variables $(N_i)_{i \geq 1}$ with common distribution ${\mathcal N} (0, \sigma^2(h))$, such that,
\[
S_n(h)- \sum_{i=1}^n N_i = o(n^{1/p}) \, \text{ ${\mathbb P}$-a.s.}
\]
\end{corollary}
\noindent {\bf Proof.}  The result comes from an application of Theorem \ref{KMTavecdeltan1} by taking into account Remark \ref{remarksufficient}. As already mentionned, 
$\nu$ admits a moment of order $S-\tau= p+(p-1)\tau$.  Hence, one can prove that condition \eqref{condondelta} holds with $r=p+\tau(p-1)$,  by using the last statement of the following lemma (taking  $\gamma =  (pr-2r+1)/(r-p) = -2 + (S-1)/\tau$).
\begin{Lemma}\label{propDR}
Let $\gamma>-1$ and $t >0$. Assume that 
$S\ge  t+(\gamma+2)\tau$. Then \[\sum_{n\ge 1} n^\gamma  \iint \E |W_{n-1,x} - W_{n-1,y}|^t \nu (d x) \nu (dy)  <\infty \, . \] In particular, for any Lipschitz function $h$,  if  $S\ge  1+(\gamma+2)\tau$ then $\sum_{n\ge 1} n^\gamma  \delta (n) <\infty$. 
\end{Lemma}
The proof of the lemma above is postponed to the Appendix (see Section \ref{sectionproofpropDR}).

\section{Proofs of Theorems \ref{KMTavecdeltaninfty} and \ref{KMTavecdeltan1}}
The proofs of Theorems \ref{KMTavecdeltaninfty} and \ref{KMTavecdeltan1} follow the scheme of proof of Theorem 2.1 in Berkes-Liu-Wu \cite{BLW14} by applying the following general Proposition \ref{generalpropBLW}, which comes from 
a careful analysis of the proof of their strong approximation result. To state this general proposition several preliminary notations are needed.

\smallskip

\noindent{\it A Preliminary result.} For  Proposition \ref{generalpropBLW}  below, we consider $(X_k)_{k \geq 1}$ a strictly stationary sequence of real-valued random variables in 
${\mathbb L}^p$ ($p>2$) and $( \varepsilon_i)_{i \geq 0}$ a sequence of iid random variables.  Let $(M_k)_{k \geq 1}$ be a sequence of positive real 
numbers and define 
\beq \label{defphik}
\varphi_k (x) = (x \wedge M_k) \vee (-M_k)  \text{ and } g_k (x) = x-\varphi_k (x)  \, .
\eeq
Then, define
\beq \label{defXkj}
{ X}_{k,j} =\varphi_k (X_j)  - \E \varphi_k (X_j) \, \mbox{ and } \, 
W_{k, \ell} = \sum_{i=1+3^{k-1}}^{\ell +  3^{k-1}} { X}_{k,i}  \, .
\eeq
Let  now $(m_k)_{k \geq 1}$ be a non-decreasing sequence of positive integers such that $m_k =o ( 3^k) $, as $k \rightarrow \infty$,  and define 
\beq \label{defXkjtilde}
{\tilde X}_{k,j} = \E \big ( \varphi_k (X_j)  | \varepsilon_j, \varepsilon_{j-1}, \ldots, \varepsilon_{j-m_k}\big ) - \E \varphi_k (X_j) \text{ for any  $j \geq m_k$} \, 
\text{ and } \, {\widetilde W}_{k, \ell} = \sum_{i=1+3^{k-1}}^{\ell +  3^{k-1}} {\tilde X}_{k,i} \, .
\eeq
Finally set $k_0 := \inf \{ k \geq 1 \, : \,  m_k \leq 2^{-1} 3^{k-2} \}$ and define
\beq \label{defnuk}
\nu_k = m_k^{-1}  \big \{ \E ({\widetilde W}^2_{k, m_k} )  + 2 \E ({\widetilde W}_{k, m_k} ({\widetilde W}_{k, 2m_k}  -{\widetilde W}_{k, m_k} ))\big \} \, .
\eeq
The general proposition coming from a careful analysis of the proof of Theorem 2.1 in Berkes-Liu- Wu \cite{BLW14} reads as follows
\begin{Proposition}[Berkes-Liu-Wu \cite{BLW14}] \label{generalpropBLW} 
Let $p>2$. Assume that we can find a sequence of positive real numbers $(M_k)_{k \geq 1}$  a  non-decreasing sequence of positive integers $(m_k)_{k \geq 1}$   such that $m_k =o ( 3^{2k/p} k^{-1}) $, as $k \rightarrow \infty$, in such a way that the following conditions are satisfied:
\beq \label{condtruncature}
\sum_{k \geq 1} 3^{k (p-1)/p}  \E (|g_k (X_1)|)< \infty \, , 
\eeq
there exists  $\alpha \geq 1$ such that 
\beq \label{condmkdependence}
\sum_{k \geq k_0} 3^{- \alpha k/p} \Big \Vert \max_{1 \leq \ell \leq 3^k - 3^{k-1}} \big |   { W}_{k, \ell} -  {\widetilde W}_{k, \ell} \big | \Big \Vert_{\alpha}^{\alpha}  < \infty \, , 
\eeq
and there exists $ r \in ]2, \infty [$ such that 
\beq \label{condtoapplySakhanenko}
\sum_{k \geq k_0} \frac { 3^{k} }{  3^{kr /p}m_k}   \E \Big (  \max_{1 \leq \ell \leq  3m_k }  \big |   {\widetilde W}_{k, \ell}    \big |^r \Big ) < \infty \, .
\eeq
Assume in addition that 
\beq \label{defsigma2asalimit}
\text{the series} \ \sigma^2 =  {\rm Var} (X_1) + 2 \sum_{i \geq 1} {\rm Cov} (X_1, X_{i+1} ) \  \text{converge,}
\eeq
and
\beq \label{cond2v_k}
3^k (  \nu_k^{1/2}  - \sigma)^2 = o(3^{2k/p} (\log k)^{-1})  \, , \, \mbox{ as $k \rightarrow \infty$}\, .  
\eeq
Then, one can redefine $(X_n)_{n \geq 1}$ without changing its distribution on a (richer) probability space on which 
there exist iid random variables $(N_i)_{i \geq 1}$ with common distribution ${\mathcal N} (0, \sigma^2)$, such that,
\beq \label{conclusiongeneralprop}
 S_n - n \E (X_1)- \sum_{i=1}^n N_i = o(n^{1/p}) \, \text{ ${\mathbb P}$-a.s.}
\eeq
\end{Proposition}
Note that  \eqref{defsigma2asalimit} implies that $n^{-1}{\rm Var} (S_n)$ converges to 
$\sigma^2$ (which is therefore non-negative). 
Let us now briefly  explain how the proposition follows from the work of Berkes-Liu-Wu \cite{BLW14}.

\medskip

Condition \eqref{condtruncature} together with condition \eqref{condmkdependence} prove that it is enough to show \eqref{conclusiongeneralprop} with 
\[
{\tilde S}_n= \sum_{k=1}^{h_n-1} {\widetilde W}_{k, 3^k - 3^{k-1}} +{\widetilde W}_{h_n, n- 3^{h_n -1}} 
\]
instead of $ S_n - n \E (X_1)$,  
where, for $n \geq 2$,  $h_n :=  \lceil (\log n)/(\log 3)\rceil $ (so that $h_n$ is the unique integer such that $3^{h_n -1} < n \leq 3^{h_n}$). Next, condition \eqref{condtoapplySakhanenko} allows first to show that the proof of the 
proposition is reduced to prove \eqref{conclusiongeneralprop} with $S_n^{\diamond} $ replacing  $S_n - n \E (X_1)$ where 
\[
S_n^{\diamond} =  \sum_{k=K_0}^{h_n-1} \sum_{j=1}^{q_k}B_{k,j}  + \sum_{j=1}^{\tau_n}B_{h_n,j}  \, , \, \text{ where } \tau_n = \Big [  \frac{ n-3^{h_n-1}}{3m_{h_n}}\Big ] -2 \, ,
\]
with
$B_{k, j} =0$ if $k < k_0$ and for $k \geq k_0$,
\[
B_{k,j} =  \sum_{i=1+3j m_k + 3^{k-1}}^{3(j+1) m_k  +  3^{k-1}} {\tilde X}_{k,i} \, , \, j=1, 2, \ldots, q_k= [ 2 \times 3^{k-2} / m_k] -2  \, .
\]
A careful analysis of the steps 3.2 and 3.3  of the proof of Theorem 2.1 in Berkes-Liu-Wu \cite{BLW14} reveals that condition \eqref{condtoapplySakhanenko}  is also sufficient to apply Theorem 1 in Sakhanenko \cite{Sa06} 
(at different steps of their proof) and this leads to the following strong approximation result: one can redefine $(X_n)_{n \geq 1}$ without changing its distribution on a (richer) probability space on which 
there exists a standard Brownian motion $B= \{ B(t), t \in {\mathbb R}^+ \}$ such that,
\[
 \max_{i \leq n} \big | S_i^{\diamond}   - B( \sigma_i^2) \big | = o(n^{1/p}) \, \text{ ${\mathbb P}$-a.s.}
\]
where
\[
\sigma_n^2 = \sum_{k=k_0}^{h_n-1} 3 m_k q_k \nu_k + 3 m_{\tau_n} \tau_n \nu_{h_n} \, .
\]
The last step 3.4 of their proof then consists  in showing that  one can construct another standard Brownian motion $W=\{ W(t), t \in {\mathbb R}^+ \}$ (depending on $B$) such that 
\[
 B( \sigma_n^2) - \sigma W(n) = o(n^{1/p}) \, \text{ ${\mathbb P}$-a.s.}
\]
This step is achieved provided that we can prove that $\nu_k \rightarrow \sigma^2$, $m_k =o ( 3^{2k/p} k^{-1}) $, as $k \rightarrow \infty$, and condition \eqref{cond2v_k} holds. 

\medskip

\noindent {\it Some preliminary considerations.} 
The following considerations allowing to extend the stationary sequence $(X_n)_{n \geq 1}$ defined by \eqref{functionofarandomiterates} to a stationary sequence on ${\mathbb Z}$ will be useful.

 \smallskip
 
 For any $n \geq 1$, let $V_n =(\varepsilon_n, W_{n-1})$. Hence $(X_n)_{n \geq 1}$ is a functional of the Markov chain $(V_n)_{n \geq 1}$ 
with state space $G \times X$ and stationary distribution $\mu \otimes \nu$.   The Markov chain $(V_n)_{n \geq 1}$ being stationary, by Kolmogorov's theorem, there exists a probability ${\hat \p}$ on the measurable space 
$({\hat \Omega}, {\hat {\mathcal F}})= ((G \times X)^{\mathbb Z}, ({\mathcal B }(G) \times{\mathcal B }( X))^{\mathbb Z} ) $ invariant by the shift ${\hat \eta}$ on ${\hat \Omega}$ and such that the law of the coordinate process 
$({\hat V}_n = ( {\hat \varepsilon}_n, {\hat W}_{n-1}) )_{n \in {\mathbb Z}}$ (with values in $G \times X$) under ${\hat \p}$ is the same as the one of $(V_n)_{n \geq 1}$ under $\p_{\nu}$.  Hence, if we define for any integer $n$, 
${\hat X}_n:= h ({\hat V}_0) \circ  {\hat \eta}^n$, it follows that $({\hat X}_n)_{n \in {\mathbb Z}}$ forms a stationary sequence with stationary distribution  $\pi$, whose law under 
${\hat \p}$ is the same as the one of $(X_n)_{n \geq 1}$ under $\p_{\nu}$. To prove the theorem, it suffices then to prove that it holds for the extended sequence  $({\hat X}_n)_{n \in {\mathbb Z}}$ which is  a 
stationary sequence adapted to the stationary filtration $({\widehat {\mathcal F}}_n)_{n \in {\mathbb Z}}$ where ${\widehat {\mathcal F}}_n= \sigma ( {\hat V}_k  , k \leq n) $
To avoid additional notations, in the rest of the proof we write $(X_n)_{n \in {\mathbb Z}}$ for $({\hat X}_n)_{n \in {\mathbb Z}}$, $(V_n)_{n \in {\mathbb Z}}$ for $({\hat V}_n)_{n \in {\mathbb Z}}$ and 
$( {\mathcal F}_n)_{n \in {\mathbb Z}}$ for $({\widehat {\mathcal F}}_n)_{n \in {\mathbb Z}}$.

\smallskip

\subsection{Proof of Theorem \ref{KMTavecdeltaninfty}} By the reverse martingale convergence theorem
and stationarity, $\Vert \E ( X_n |  {\mathcal F}_0)  - \E (X_n)\Vert_2$ is decreasing to $\Vert \E ( X_0 |  {\mathcal F}_{- \infty})  - \E (X_0)\Vert_2$, as $n \rightarrow \infty$. Hence, by condition \eqref{condondeltaninfty}, 
$\E ( X_0 |  {\mathcal F}_{- \infty}) =  \E (X_0)$ a.s. Applying Lemma \ref{lmapourcov} of the Appendix and taking into account condition \eqref{condondeltaninfty}, we get (since $q>1/2$),  
\[
\sum_{k \geq 1}  |   {\rm Cov} (X_1, X_{k+1} )| \ll  \Big (  \sum_{k \geq 0}  (k+1)^{-1/2} \Vert   \E ( X_{k}  | V_0 )  - \E (X_k) \Vert_2  \Big )^2< \infty \, . 
\]
This proves that  the series $\sigma^2 =  {\rm Var} (X_1) + 2 \sum_{i \geq 1} {\rm Cov} (X_1, X_{i+1} )$ converge
absolutely and  condition 
\eqref{defsigma2asalimit}  of Proposition \ref{generalpropBLW} holds.

\medskip

Assume first that $\sigma^2 =0$. To prove that $S_n - n \E (X_1) =o(n^{1/p})$ a.s., we shall use Theorem 4.7 in Cuny-Merlev\`ede \cite{CM14}. Hence, it suffices to prove that 
\beq \label{toapplyTh47}
\sum_{n \geq 1} \frac{\Vert S_n  - n \E (X_1)\Vert_p}{n^{1 + 1/p}} < \infty \, .
\eeq
With this aim, we start by noticing that by condition \eqref{condondeltaninfty},
\[
\sum_{k \geq 1} k^{-1/p} \Vert \E ( X_k | V_0) - \E (X_k)\Vert_p < \infty   \, \mbox{ and therefore } \, \sum_{k \geq 1}  \frac{ \Vert \E ( S_k | V_0) -  k \E (X_1)\Vert_p}{k^{1 + 1/p}} < \infty \, .
\]
Theorem 2.3 in  \cite{CM14} then asserts that there exists a stationary sequence $(D_k)_{k \in {\mathbb Z}}$ of martingale differences in ${\mathbb L}^p$, adapted to 
$({\mathcal F}_k)_{k \in {\mathbb Z}}$ and such that $n^{-1/2} \Vert S_n - n \E (X_1) - \sum_{k=1}^n D_k \Vert_p \rightarrow 0$, as $n \rightarrow \infty$. Together with the fact that 
 $\lim_{n \rightarrow \infty} n^{-1}{\rm Var} (S_n) = \sigma^2 =0$, it follows that $D_k =0$ a.s, for any $k$. Therefore,  the upper bound (4) in  \cite{CM14} and condition \eqref{condondeltaninfty} entail that
\begin{multline*}
\sum_{n \geq 1} \frac{\Vert S_n - n \E (X_1)\Vert_p}{n^{1 + 1/p}}  \ll \sum_{n \geq 1} \frac{n^{1/2} }{n^{1 + 1/p}}  \sum_{k \geq [n^{p/2}]}  \frac{ \Vert \E ( S_k | V_0) -  k \E (X_1)\Vert_p}{k^{1 +1/p}}  \\
\ll  \sum_{k \geq 1} \frac{1 }{k^{1 + 2/p^2}}  \Vert \E ( S_k | V_0) -  k \E (X_1)\Vert_p  
\ll  \sum_{k \geq 1} \frac{1 }{k^{ 2/p^2}}  \Vert \E ( X_k | V_0) -   \E (X_1)\Vert_p  \ll \sum_{k \geq 1} \frac{1 }{k^{ q+ 2/p^2}}   \, ,
\end{multline*}
which proves \eqref{toapplyTh47} since $ q+ 2/p^2 -1 > (2 p^2)^{-1} (p^3 - 3 p^2 +4 ) = (2 p^2)^{-1} (p-2)^2 (p+1) >0$. The theorem is then proved in the case where $\sigma^2 =0$. 

\medskip

Assume from now that $\sigma^2 >0$. We choose
\[
M_k =3^{k/p} \text{ and } m_k = [3^{ 2 ( 1 - \varepsilon) k/p}  ]  \text{ with } 0 < \varepsilon < \min \Big ( 1 - \frac{p-1}{ 2q} , \frac{1}{2} \Big )\, ,
\]
Note that the sequence $(m_k)_{k \geq 0 }$ satisfies $m_k =o ( 3^{2k/p} k^{-1}) $, as $k \rightarrow \infty$.  We prove below that conditions \eqref{condtruncature}, \eqref{condmkdependence}, \eqref{condtoapplySakhanenko} and \eqref{cond2v_k}
of Proposition \ref{generalpropBLW} are satisfied with the above choices of $(M_k)_{k \geq 0}$ and $(m_k)_{k \geq 0}$. 

\medskip

Since the $X_i$'s are in ${\mathbb L}^p$, it is easy to see that with the choice of $M_k$, condition  \eqref{condtruncature} is satisfied (it suffices to write that  $\E (|g_k (X_1)|)\leq   \E (|X_1|{\mathbf 1}_{|X_1| >M_k})$ and to use Fubini's Theorem).  
Next,  for $k \geq k_0$,  
Lemma \ref{lmaenlevertildedeltan} of the Appendix combined with condition \eqref{condondeltaninfty} implies that
\[
\Big \Vert  \max_{1 \leq \ell \leq 3^k - 3^{k-1}} \big |   { W}_{k, \ell} -  {\widetilde W}_{k, \ell} \big | \Big \Vert_1 \leq  \sum_{i=1+3^{k-1}}^{ 3^{k}} \Vert X_{k,i}  - {\tilde X}_{k,i}  \Vert_1 
\leq  \frac{C 3^k }{m_k^{q}}  \, .
\]
Therefore,
\[
\sum_{k \geq k_0} 3^{-k/p} \Big \Vert  \max_{1 \leq \ell \leq 3^k - 3^{k-1}} \Big |   { W}_{k, \ell} -  {\widetilde W}_{k, \ell} \big | \big \Vert_1 \ll   \sum_{k \geq 1} \frac{3^{k (p-1)/p}}{ 3^{ 2 q ( 1 - \varepsilon) k/p}}  < \infty  \, ,
\]
since $ 2 q ( 1 - \varepsilon)  >p-1 $. Condition \eqref{condmkdependence} is then satisfied with $\alpha =1$. We prove now that  we can find a real 
number $ r \in ]2, \infty [$ such that 
 \eqref{condtoapplySakhanenko} holds. Let $r \geq 2$, 
\[
Y_{k,i} = {\tilde X}_{k, i+3^{k-1}} \, , \, {\mathcal G}_{k,i} := {\mathcal G}_{i + 3^{k-1} } = \sigma ( \varepsilon_{i + 3^{k-1}}, \ldots,  \varepsilon_{ 1 }, W_0 )  \, , 
\]
and
\[
{d_{k,i}} = Y_{k,i}  - \E ( Y_{k,i}  |  {\mathcal G}_{k,i-1} ) \, .
\]
With these notations, we have
\beq  \label{dec1tildeWnormr}
 \Big \Vert  \max_{1 \leq \ell \leq  3m_k }  \big |   {\widetilde W}_{k, \ell}    \big | \Big \Vert_r \leq  \Big \Vert  \max_{1 \leq \ell \leq  3m_k }  \Big |  \sum_{i=1}^{\ell} {d_{k,i}}    \Big | \Big \Vert_r + 
 \Big \Vert  \max_{1 \leq \ell \leq  3m_k }  \Big |  \sum_{i=1}^{\ell}  \E ( Y_{k,i}  |  {\mathcal G}_{k,i-1} )    \Big | \Big \Vert_r  \, .
\eeq
By Rosenthal's inequality for martingales, 
\[
 \Big \Vert  \max_{1 \leq \ell \leq  3m_k }  \Big |  \sum_{i=1}^{\ell} {d_{k,i}}    \Big | \Big \Vert_r^r \ll \Big ( \sum_{i=1}^{3m_k} \Vert \E ( { d^2_{k,i}} |  {\mathcal G}_{k,i-1} )  \Vert_{r/2}  \Big )^{r/2} + 
  \sum_{i=1}^{3m_k} \Vert d_{k,i}\Vert^r_{r} \, .
\]
Note that 
\[
 \E ( { d^2_{k,i}} |  {\mathcal G}_{k,i-1} )  \leq  \E ( \E^2 ( \varphi_k (X_{i+3^{k-1}} )  | {\mathcal H}_{k,i})  |  {\mathcal G}_{k,i-1} ) \leq  \E ( \E ( X^2_{i+3^{k-1}}  | {\mathcal H}_{k,i})  |  {\mathcal G}_{k,i-1} ) \, , 
\]
where ${\mathcal H}_{k,i} =  \sigma ( \varepsilon_{i + 3^{k-1}}, \ldots,  \varepsilon_{  i + 3^{k-1} -m_k} )  $.  Here, recall  the following well known fact: 
if $Y$ is an integrable random variable, and ${\mathcal G}_1$ and ${\mathcal G}_2$ are two  $\sigma$-algebras such that $ \sigma(Y) \vee {\mathcal G}_1$ is independent of  ${\mathcal G}_2$, then
\beq \label{lmatribu}
\E ( Y | {\mathcal G}_1 \vee {\mathcal G}_2)  =\E ( Y | {\mathcal G}_1) \ \text{ a.s.}
\eeq
Applying \eqref{lmatribu} with 
${\mathcal G}_1=\sigma ( \varepsilon_{i + 3^{k-1}-1}, \ldots,  \varepsilon_{  i + 3^{k-1} -m_k} ) $, ${\mathcal G}_2= {\mathcal G}_{k,i -m_k-1}  $ and $Y =  \E ( X^2_{i+3^{k-1}}  | {\mathcal H}_{k,i}) $, 
we get 
\[
 \E ( { d^2_{k,i}} |  {\mathcal G}_{k,i-1} )  \leq   \E ( X^2_{i+3^{k-1}}  | \sigma ( \varepsilon_{i + 3^{k-1}-1}, \ldots,  \varepsilon_{  i + 3^{k-1} -m_k} ) )  \, .
\]
Hence, by assumption \eqref{condonnorminftycarre}, 
\[
\Vert \E (  d^2_{k,i}  |  {\mathcal G}_{k,i-1} )  \Vert_{r/2}  \leq \Vert \E ( X^2_{i+3^{k-1}} |  {\mathcal G}_{k,i-1} )  \Vert_{r/2} \leq \Vert \E ( X^2_{i+3^{k-1}} |  {\mathcal G}_{k,i-1} )  \Vert_{\infty} \leq c \, .
\]
On another hand, by stationarity, 
\[
\Vert d_{k,i}\Vert^r_{r} \ll \Vert  \varphi_k (X_{0}) \Vert^r_{r} \ll \Vert X_{0} {\mathbf 1}_{|X_0| \leq M_k } \Vert^r_{r} + M^r_k \p (|X_0| > M_k) \, .
\]
So, overall, 
\beq  \label{dec1tildeWnormrterm1}
  \Big \Vert  \max_{1 \leq \ell \leq  3m_k }  \Big |  \sum_{i=1}^{\ell} {d_{k,i}}    \Big | \Big \Vert^r_r \ll m_k^{r/2} + m_k \big (  \Vert X_{0} {\mathbf 1}_{|X_0| \leq M_k } \Vert^r_{r} + M^r_k \p (|X_0| > M_k) \big ) \, .
 \eeq
We handle now the second term in the right-hand side of \eqref{dec1tildeWnormr}. We apply Proposition \ref{LpmainNS}  of the Appendix with $\alpha =r$, $r=r_k$ where 
$r_k$ is the unique positive integer such that $2^{r_k-1} \leq 3 m_k < 2^{r_k}$, 
\[
Z_i = 0 \mbox{ for $i \leq 0$ , } \, Z_{i}:=Z_{k,i} = \E ( Y_{k,i}  |  {\mathcal G}_{k,i-1} )   \mbox{ for $i \geq 1 $}  
\]
and
\[
 {\mathcal F}_{i}= \{\emptyset, \Omega\} \mbox{ for $i \leq 0$ , } \,   {\mathcal F}_{i} = \sigma ( \varepsilon_{i  -1 + 3^{k-1} }, \ldots , \varepsilon_1, W_0 ) = {\mathcal G}_{k,i-1} \mbox{ for $i \geq 1 $}\, .
\]
We then get 
\begin{multline} \label{appli1propLpMainNS}
 \Big \Vert  \max_{1 \leq \ell \leq  3m_k }  \Big |  \sum_{i=1}^{\ell}  \E ( Y_{k,i}  |  {\mathcal G}_{k,i-1} )    \Big | \Big \Vert_r \\ \ll  \Big (  \sum_{j=1}^{3m_k} \Vert  \E ( Y_{k,j}  |  {\mathcal G}_{k,j-1} )   \Vert_r^2 \Big )^{1/2} 
 + 
 \sum_{\ell=0}^{r_k-1}  \Big ( \sum_{m=1}^{2^{r_k-\ell}} \Vert {\mathbb E} ( T_{m 2^\ell } -T_{(m-1) 2^\ell }  | {\mathcal F}_{ (m-2) 2^\ell +1} ) \Vert_r^2 \Big )^{1/2} \, ,
\end{multline}
where $T_{\ell} = \sum_{i=1}^{\ell}   \E ( Y_{k,i}  |  {\mathcal G}_{k,i-1} )$. 
By fact \eqref{lmatribu}, we note that, for any $i \geq 1$, 
\begin{multline*}
 \Vert\E ( Y_{k,i}  |  {\mathcal G}_{k,i-1} )  \Vert_{\infty} = \Vert   \E ( X_{k,i+3^{k-1}}  | \sigma ( \varepsilon_{i + 3^{k-1}-1}, \ldots,  \varepsilon_{  i + 3^{k-1} -m_k} ) ) \Vert_{\infty}  \\
 \leq 2  \Vert\E ( \varphi_k (X_{i + 3^{k-1}} ) |  {\mathcal G}_{k,i-1} )  \Vert_{\infty} \leq 2  \Vert\E ( \varphi^2_k (X_{i + 3^{k-1}} ) |  {\mathcal G}_{k,i-1} )  \Vert_{\infty}^{1/2}  \, .
\end{multline*}
Therefore, by condition \eqref{condonnorminftycarre}, 
\beq \label{appli1propLpMainNSb1}
 \Vert\E ( Y_{k,i}  |  {\mathcal G}_{k,i-1} )  \Vert_{\infty} \leq 2 \sqrt{c} \, .
\eeq
Next, since  ${\mathcal F}_{i}= \{\emptyset, \Omega\}$  for $i \leq 0$ and the $Z_i$'s are centered ,  for any $\ell \geq 0$, 
\[
\Vert {\mathbb E} ( T_{2^{\ell} }   | {\mathcal F}_{-2^{\ell}+1} ) \Vert_r =0 \, . 
\]
Moreover, for any $m \geq 2$ and any $\ell \geq 0$, 
\[
\Vert {\mathbb E} ( T_{m 2^\ell } -T_{(m-1) 2^\ell }  | {\mathcal F}_{ (m-2) 2^\ell +1 } ) \Vert_r \leq \sum_{i=(m-1)2^{\ell}+1}^{m 2^{\ell}} \Vert  \E (   {\tilde X}_{k,i+3^{k-1}} | {\mathcal G}_{ k, (m-2) 2^\ell   } )\Vert_r \, .
\]
But, for any $m \geq 2$, any $\ell \geq 0$  and any $i \geq (m-1)2^{\ell}+1$, 
\begin{multline*}
\Vert  \E (  {\tilde X}_{k,i+3^{k-1}} | {\mathcal G}_{ k, (m-2) 2^\ell   } ) \Vert_r \\ =
 \Vert  \E ( \E ( \varphi_k (X_{i+3^{k-1}} ) | \varepsilon_{i+3^{k-1}}, \ldots,  \varepsilon_{i+3^{k-1} -m_k}) | {\mathcal G}_{ k, (m-2) 2^\ell   } ) - \E (\varphi_k (X_{i+3^{k-1}})) \Vert_r \, .
\end{multline*}
Hence, if $2^{\ell } \geq m_k$, 
\[
\Vert  \E (  {\tilde X}_{k,i+3^{k-1}} | {\mathcal G}_{ k, (m-2) 2^\ell   } ) \Vert_r = 0 \, ,
\]
and if $2^{\ell } \leq  m_k-1$,  by using \eqref{lmatribu}, 
\[
\Vert  \E (  {\tilde X}_{k,i+3^{k-1}} | {\mathcal G}_{ k, (m-2) 2^\ell   } ) \Vert_r  \leq \Vert  \E (  { X}_{k,i+3^{k-1}} | {\mathcal G}_{ k, (m-2) 2^\ell   } ) \Vert_{\infty}  \, .
\]
But, by using stationarity,  the Markov property and the fact that $\varphi_k$ is 1-Lipschitz, 
\begin{multline*}
\Vert  \E (  { X}_{k,i+3^{k-1}} | {\mathcal G}_{ k, (m-2) 2^\ell   } ) \Vert_{\infty} \leq  \sup_{x \in X} \Big |  \E ( \varphi_k (X_{i- ( (m-2) 2^\ell , x })  )   - \int   \E ( \varphi_k (X_{i- ( (m-2) 2^\ell , y })  )  d \nu (y) \Big | \\
\leq  \sup_{x \in X}  \int  \big |  \E ( \varphi_k (X_{i- ( (m-2) 2^\ell , x })  )   -   \E ( \varphi_k (X_{i- ( (m-2) 2^\ell , y })  ) \big |  d \nu (y) \, .
\end{multline*}
Hence, for any $m \geq 2$, any $\ell \geq 0$  and any $i \geq (m-1)2^{\ell}+1$,
\beq \label{majinftyphik}
\Vert  \E (  {\tilde X}_{k,i+3^{k-1}} | {\mathcal G}_{ k, (m-2) 2^\ell   } ) \Vert_r \leq  \Vert  \E (  { X}_{k,i+3^{k-1}} | {\mathcal G}_{ k, (m-2) 2^\ell   } ) \Vert_{\infty}   \ll \frac{1}{ (  i - (m-2) 2^\ell )^q } \, .
\eeq
Since $q>1/2$, the above considerations imply that 
\[
 \sum_{\ell=0}^{r_k-1}  \Big ( \sum_{m=1}^{2^{r_k-\ell}} \Vert {\mathbb E} ( T_{m 2^\ell } -T_{(m-1) 2^\ell }  | {\mathcal F}_{(m-2) 2^\ell +1} ) \Vert_r^2 \Big )^{1/2}  \ll 
2^{r_k/2}  \sum_{\ell=0}^{r_k-1}  2^{ \ell /2} 2^{-\ell q  }  \ll m_k^{1/2}\, .
\]
Combined with \eqref{appli1propLpMainNS} and \eqref{appli1propLpMainNSb1}, the upper bound above implies that
\beq  \label{dec1tildeWnormrterm2}
 \Big \Vert  \max_{1 \leq \ell \leq  3m_k }  \Big |  \sum_{i=1}^{\ell}  \E ( Y_{k,i}  |  {\mathcal G}_{k,i-1} )    \Big | \Big \Vert^r_r \ll m_k^{r/2} \, .
 \eeq
Hence, starting from \eqref{dec1tildeWnormr} and taking into account \eqref{dec1tildeWnormrterm1} and \eqref{dec1tildeWnormrterm2},  we get that for any $r \geq 2$, 
\[
 \Big \Vert  \max_{1 \leq \ell \leq  3m_k }  \Big |   {\widetilde W}_{k, \ell}    \big | \big \Vert^r_r \ll m_k^{r/2} +m_k \big (   \Vert X_{0} {\mathbf 1}_{|X_0| \leq M_k } \Vert^r_{r} + M^r_k \p (|X_0| > M_k) \big )  \, .
\]
This implies that \eqref{condtoapplySakhanenko} holds with $r > \max \big \{ 2,  \varepsilon^{-1} \big ( p-2 (1 - \varepsilon) \big ) \big \}$. 
\smallskip

To end the proof it remains to prove condition \eqref{cond2v_k}. Note first that since $\sigma^2$ is assumed to be strictly positive, we have
\[
 | \nu_k^{1/2}  - \sigma |= \sigma^{-1} \sigma | \nu_k^{1/2}  - \sigma | \leq \sigma^{-1}  (\nu_k^{1/2}  + \sigma )  | \nu_k^{1/2}  - \sigma |  =
 \sigma^{-1}    | \nu_k  - \sigma^2 | \, ,   \]
and therefore condition \eqref{cond2v_k} reads as 
\beq \label{cond2v_kbis}
3^k (  \nu_k - \sigma^2)^2 = o(3^{2k/p} (\log k)^{-1})  \, , \, \mbox{ as $k \rightarrow \infty$}\, .  
\eeq
To verify condition \eqref{cond2v_kbis}, let us define, for $i \geq 0$, 
\[
{\tilde c}_{k,i}= {\rm cov} ({\tilde X}_{k, m_k+1} , {\tilde X}_{k,i+m_k+1} ) \, \text{ and } {\hat c}_{k,i} = {\rm cov} ({ X}_{k,0} , { X}_{k,i} ) \, .
\]
Using stationarity, we have
\[
 \E ({\widetilde W}^2_{k, m_k} )  = m_k {\tilde c}_{k,0} + 2 \sum_{\ell =1}^{m_k -1} (m_k - \ell) {\tilde c}_{k,\ell} \,  \text{ and } \,  \E ({\widetilde W}_{k, m_k} ({\widetilde W}_{k, 2m_k}  -{\widetilde W}_{k, m_k} )) =  \sum_{\ell =1}^{m_k}  \ell   {\tilde c}_{k,\ell}\, .
\]
Therefore
 \beq \label{nuk1}
\nu_k = {\tilde c}_{k,0}  + 2 \sum_{\ell =1}^{m_k }  {\tilde c}_{k,\ell} \, . 
\eeq
We first prove that 
\beq \label{compnukavechat}
\Big | \nu_k -  \big (  {\hat c}_{k,0}  + 2 \sum_{\ell \geq 1 }  {\hat c}_{k,\ell}  \big )  \Big |  \ll  3^{k/(2p)} m_k^{-q/2} + \big (  m_k^{-1/2} ( \log m_k){\bf 1}_{q =1}  + m_k^{-q +1/2} {\bf 1}_{q <1} \big )  \, .
\eeq
With this aim we use the arguments developed in \cite{BLW14}  to get their inequality (3.56). Hence, we start by noting that since $\varphi_k$ is $1$-Lipschitz, 
$ \big ( \Vert \E (\varphi_k (X_n) | {\mathcal F}_0) - \E ( \varphi_k (X_n) \Vert_2) \big  )_{n \geq 0}$ is a decreasing sequence such that  $ \Vert \E (\varphi_k (X_n) | {\mathcal F}_0) - \E ( \varphi_k (X_n) \Vert_2 \leq \delta_{\infty} (n)$. Hence, by the same arguments as those developed in the first lines of the proof of Theorem \ref{KMTavecdeltaninfty}, we infer that, under condition \eqref{condondeltaninfty}, there exists a constant $C$ not depending on $k$ such that $\sum_{\ell \in {\mathbb Z} }  | {\hat c}_{k,\ell} | \leq C$. 
Therefore,  $\lim_{j \rightarrow \infty} j^{-1} \E ({ W}^2_{k, j}) 
=  {\hat c}_{k,0}  + 2 \sum_{\ell \geq 1 }  {\hat c}_{k,\ell} $. On another hand,  the following convergence clearly holds:  $\lim_{j \rightarrow \infty} j^{-1} \E ({\widetilde W}^2_{k, j}) = \nu_k $. In addition,   for all $j \geq 1$, 
\[
|\E ({\widetilde W}^2_{k, j}) - \E ({ W}^2_{k, j}) | \leq \Vert {\widetilde W}_{k, j} - W_{k, j}\Vert_2 \Vert {\widetilde W}_{k, j} + W_{k, j}\Vert_2 \, .
\]
The above considerations imply 
\begin{equation} \label{dec1compnukavechat}
\Big | \nu_k -  \big (  {\hat c}_{k,0}  + 2 \sum_{\ell \geq 1 }  {\hat c}_{k,\ell}  \big )  \Big |  \leq 
 \limsup_{j \rightarrow \infty} j^{-1} \Vert {\widetilde W}_{k, j} - W_{k, j}\Vert_2^2  + 2 \sqrt{ C}  \limsup_{j \rightarrow \infty} j^{-1/2} \Vert {\widetilde W}_{k, j} - W_{k, j}\Vert_2 \, .
\end{equation}
To take care of $\Vert {\widetilde W}_{k, j} - W_{k, j}\Vert_2$, we apply  Proposition \ref{LpmainNS} of the Appendix with, this time,  $\alpha =2$,  $r=r_j$ where $r_j$ is the unique positive integer such that $2^{r_j-1} \leq j < 2^{r_j}$, 
\[
Z_i = 0 \mbox{ for $i \leq 0$ , } \, Z_{i}:=Z_{i,k} = X_{k,i+3^{k-1}} - {\tilde X}_{k, i+3^{k-1}}  \mbox{ for $i \geq 1 $} \, ,   
\]
and
\[
 {\mathcal F}_{i}= \{\emptyset, \Omega\} \mbox{ for $i \leq 0$ , } \,   {\mathcal F}_{i} = \sigma ( \varepsilon_{i + 3^{k-1} }, \ldots , \varepsilon_1, W_0 ) = {\mathcal G}_{  k,i }  \mbox{ for $i \geq 1 $}\, .
\]
Hence
\begin{multline} \label{appli2propLpMainNS}
\Vert {\widetilde W}_{k, j} - W_{k, j}\Vert_2  \ll  \Big (  \sum_{i=1}^{j} \Vert  X_{k,i+3^{k-1}} - {\tilde X}_{k, i+3^{k-1}}  \Vert_2^2 \Big )^{1/2}  \\
 + 
 \sum_{\ell=0}^{r_k-1}  \Big ( \sum_{m=1}^{2^{r_k-\ell}} \Vert {\mathbb E} ( T_{m 2^\ell } -T_{(m-1) 2^\ell }  | {\mathcal F}_{ (m-2) 2^\ell +1} ) \Vert_2^2 \Big )^{1/2} \, ,
\end{multline}
where $T_{\ell} =  \sum_{i=1}^{\ell} (X_{k,i+3^{k-1}} - {\tilde X}_{k, i+3^{k-1}} )$. Lemma \ref{lmaenlevertildedeltan} of the Appendix  combined with condition \eqref{condondeltaninfty} implies that
\beq  \label{appli2propLpMainNSb1}
 \sum_{i=1}^{j} \Vert  X_{k,i+3^{k-1}} - {\tilde X}_{k, i+3^{k-1}}  \Vert_2^2 \leq 2 \times 3^{k/p} \sum_{i=1+3^{k-1}}^{ j+3^{k-1}} \Vert X_{k,i}  - {\tilde X}_{k,i}  \Vert_1 
\leq  \frac{C j 3^{k/p} }{m_k^{q}}  \, .
\eeq
Next, since  ${\mathcal F}_{i}= \{\emptyset, \Omega\}$  for $i \leq 0$,  for any $\ell \geq 0$, 
\[
\Vert {\mathbb E} ( T_{2^{\ell} }   | {\mathcal F}_{-2^{\ell}+1} ) \Vert_2 = | {\mathbb E} ( T_{2^{\ell} } )| = 0 \, . 
\]
Moreover, by \eqref{majinftyphik}, we infer that for any $m \geq 2$, any $\ell \geq 0$ and any $i \geq (m-1)2^{\ell}+1$, 
\begin{multline*}
\Vert  \E ( X_{k,i+3^{k-1}} - {\tilde X}_{k, i+3^{k-1}}  | {\mathcal G}_{ k, (m-2) 2^\ell +1  } )\Vert_2  \\ \leq \Vert  \E ( X_{k,i+3^{k-1}} - {\tilde X}_{k, i+3^{k-1}}  | {\mathcal G}_{ k, (m-2) 2^\ell +1  } )\Vert_{\infty} \ll \frac{1}{ (  i - (m-2) 2^\ell )^q } \ll 2^{-\ell q }\, .
\end{multline*}
On another  hand, Lemma \ref{lmaenlevertildedeltan} of the Appendix combined with condition \eqref{condondeltaninfty} implies that 
\[
\Vert  \E ( X_{k,i+3^{k-1}} - {\tilde X}_{k, i+3^{k-1}}  | {\mathcal G}_{ k, (m-2) 2^\ell +1  } )\Vert_1 \leq \Vert X_{k,i+3^{k-1}} - {\tilde X}_{k, i+3^{k-1}}  \Vert_1 \ll m_k^{-q} \, .
\]
Hence, for any $m \geq 2$, any $\ell \geq 0$ and any $i \geq (m-1)2^{\ell}+1$, we also have
\[
\Vert  \E ( X_{k,i+3^{k-1}} - {\tilde X}_{k, i+3^{k-1}}  | {\mathcal G}_{ k, (m-2) 2^\ell +1  } )\Vert^2_2   \leq 
2^{-\ell q } m_k^{-q} \, .
\]
The considerations above imply that, for any $m \geq 2$, any $\ell \geq 0$ and any $i \geq (m-1)2^{\ell}+1$, 
\[
\Vert  \E (X_{k,i+3^{k-1}} - {\tilde X}_{k, i+3^{k-1}}  | {\mathcal G}_{ k, (m-2) 2^\ell +1  } )\Vert_2 \leq 2^{-\ell q /2 } \min ( 2^{-\ell q /2 } , m_k^{-q/2}  )  \, .
\]
Hence, since $q>1/2$, 
\begin{multline} \label{appli2propLpMainNSb2}
 \sum_{\ell=0}^{r_j-1}  \Big ( \sum_{m=1}^{2^{r_j-\ell}} \Vert {\mathbb E} ( T_{m 2^\ell } -T_{(m-1) 2^\ell }  | {\mathcal F}_{(m-2) 2^\ell +1} ) \Vert_2^2 \Big )^{1/2}  \ll 
2^{r_j/2}  \sum_{\ell=0}^{r_j-1}  2^{ \ell /2} 2^{-\ell q /2 } \min ( 2^{-\ell q /2 } , m_k^{-q/2}  )  \\ \ll j^{1/2}  \big ( m_k^{-q/2} {\bf 1}_{q >1} + m_k^{-1/2} (\log m_k)  {\bf 1}_{q =1}  + m_k^{-q +1/2} {\bf 1}_{q <1} \big ) \, .
\end{multline}
Starting from \eqref{appli2propLpMainNS} and considering the upper bounds \eqref{appli2propLpMainNSb1} and \eqref{appli2propLpMainNSb2}, we get
\beq \label{step0compnukavechat}
j^{-1/2} \Vert {\widetilde W}_{k, j} - W_{k, j}\Vert_2 \ll  3^{k/(2p)} m_k^{-q/2} + \big (  m_k^{-1/2} ( \log m_k ) {\bf 1}_{q =1}  + m_k^{-q +1/2} {\bf 1}_{q <1} \big ) \, .
\eeq
Hence starting from \eqref{dec1compnukavechat} and taking into account  \eqref{step0compnukavechat} together with the fact that $3^{k/(2p)} m_k^{-q/2} \leq 2^{q/2} $, the upper bound \eqref{compnukavechat} follows. 

Let now $c_i = {\rm Cov} (X_0,X_i) $ and note that  (see Relation (3.54) in \cite{BLW14}, where the same truncation level is used)
\[
\sup_{i  \geq 0} | {\hat c}_{k,i} -  c_i  | = o (3^{-k (p-2)/p}) \, .
\]
Let 
\[
\ell_k = [ 3^{k (p-2)/ (2p)} (\log k)^{-1/2} ] \, .
\]
Since $\sigma^2 =c_0 + 2 \sum_{i \geq 1} c_i$, it follows that 
\beq \label{comp1sigma2avechat}
\Big | \sigma^2 -  (  {\hat c}_{k,0} + 2  \sum_{i \geq 1}   {\hat c}_{k,i} )  \Big | \leq o ( \ell_k 3^{-k (p-2)/p}) + 2 \sum_{i > \ell_k} | c_{i} - { \hat c}_{k,i}| \, .
\eeq
But
\[
| c_{i} - { \hat c}_{k,i}|  = |\cov ( X_0 - \varphi_k (X_0), X_i) + \cov (  \varphi_k (X_0), X_i - \varphi_k (X_i) )  | \, .
\]
Set $g_k (x) = x - \varphi_k(x)$ and note that, by the reverse martingale convergence theorem and condition  \eqref{condondeltaninfty}, $ \E (g_k (X_0) | {\mathcal F}_{- \infty} )=\E (g_k (X_0) ) $ a.s. and $ \E (X_0 | {\mathcal F}_{- \infty}) =\E (X_0) $ a.s. Hence, applying Lemma \ref{lmapourcov} of the Appendix and taking into account condition \eqref{condondeltaninfty}, we get 
\begin{multline*}
 \sum_{i > \ell_k} |\cov ( X_0 - \varphi_k (X_0), X_i)   | \ll  \sum_{j=0}^{\infty} \Vert P_{0} ( g_k (X_j)) \Vert_2 \sum_{i \geq [2^{-1} (\ell_k+ j)] +1 } 
i^{-1/2} \Vert \E (X_i | V_0) \Vert_2 \\ 
 \ll \sum_{j=0}^{\infty} \Vert P_{0} ( g_k (X_j)) \Vert_2  (\ell_k +j )^{1/2 - q } \, ,
\end{multline*}
where $P_0 (\cdot) = \E ( \cdot | {\mathcal F}_0 ) -  \E ( \cdot | {\mathcal F}_{-1} )$.  
But, by Lemma \ref{lmapourcov} of the Appendix, 
\begin{multline*}
 \sum_{j=0}^{\infty} \Vert P_{0} ( g_k (X_j)) \Vert_2  (\ell_k +j )^{1/2 - q } \ll \ell_k^{1/2 - q }  \sum_{j=0}^{\ell_k} \Vert P_{0} ( g_k (X_j)) \Vert_2  + \sum_{j \geq \ell_k} j^{1/2 - q } 
 \Vert P_{0} ( g_k (X_j)) \Vert_2 \\
 \ll \ell_k^{1/2 - q }  \sum_{j=0}^{\ell_k} (j+1)^{-1/2} \Vert  \E (  g_k (X_j) | V_0) -  \E (  g_k (X_j))  \Vert_2 + \sum_{j \geq \ell_k} j^{1/2 - q } 
 \Vert P_{0} ( g_k (X_j)) \Vert_2\, .
\end{multline*}
Note now that, since $q >1/2$,
\begin{multline}  \label{aboundforP0}
\sum_{j \geq \ell_k} j^{1/2 - q }  \Vert P_{0} ( g_k (X_j)) \Vert_2  = \sum_{j \geq \ell_k} j^{1/2 - q }  \Vert P_{0} ( g_k (X_j)) \Vert_2 \Big (  j^{-1} \sum_{m =1}^j 1  \Big )\\ = \sum_{m \geq 1}   \sum_{j \geq  \max( m , \ell_k)} j^{ - (1/2 + q) }  \Vert P_{0} ( g_k (X_j)) \Vert_2   \\
\ll   \ell_k^{1 - q }   \Big (  \sum_{j \geq   \ell_k } \Vert P_{0} ( g_k (X_j)) \Vert^2_2  \Big )^{1/2} + \sum_{m >\ell_k } m^{-q}  \Big (  \sum_{j \geq  m } \Vert P_{0} ( g_k (X_j)) \Vert^2_2  \Big )^{1/2}  \\
 \ll \ell_k^{1- q }   \Vert  \E (  g_k (X_{\ell_k}) | V_0) -  \E (  g_k (X_{\ell_k}))  \Vert_2 +  \sum_{m >\ell_k } m^{-q}   \Vert  \E (  g_k (X_m) | V_0) -  \E (  g_k (X_m))  \Vert_2  \, .
\end{multline}
So, overall, 
\begin{multline*}
 \sum_{i > \ell_k} |\cov ( X_0 - \varphi_k (X_0), X_i)   | \ll  \ell_k^{1/2 - q }  \sum_{j=0}^{\ell_k} (j+1)^{-1/2} \Vert  \E (  g_k (X_j) | V_0) -  \E (  g_k (X_j))  \Vert_2  \\ + 
\ell_k^{1- q }   \Vert  \E (  g_k (X_{\ell_k}) | V_0) -  \E (  g_k (X_{\ell_k}))  \Vert_2  
+  \sum_{m >\ell_k } m^{-q}   \Vert  \E (  g_k (X_m) | V_0) -  \E (  g_k (X_m))  \Vert_2    \, .
\end{multline*} 
Next, we note that 
\[
\Vert  \E (  g_k (X_0) | V_0) -  \E (  g_k (X_0))  \Vert_2 \ll o(3^{- k (p-2)/(2p)}) \, ,
\]
and that, for $j \geq 1$,  by condition  \eqref{condondeltaninfty}, 
\beq \label{boundnorm2gk}
\Vert  \E (  g_k (X_j) | V_0) -  \E (  g_k (X_j))  \Vert_2 \ll \min \big ( j^{-q} , j^{-q/2} 3^{- k (p-1)/(2p)} \big ) \, .
\eeq
Hence, since $q>1/2$, we infer that 
\begin{multline} \label{comp2sigma2avechat}
 \sum_{i > \ell_k} |\cov ( X_0 - \varphi_k (X_0), X_i)   |
  \ll  \ell_k^{1/2 - q }   3^{- k (p-2)/(2p)}  +  \ell_k^{1-  2q }  {\bf 1}_{q> 1} +  k \,  \ell_k^{- 1/2  }   3^{- k (p-1)/(2p)} {\bf 1}_{q=1}  
  \\
 + \big ( k \,  3^{- k (p-1)/(2p)} + \ell_k^{ 1- 3 q/2  }   3^{- k (p-1)/(2p)} + 3^{k (1 -2 q) (p-1)/ ( p q )} \big )  {\bf 1}_{q< 1}   \, .
\end{multline}
We handle now the series
\[
 \sum_{i > \ell_k} |\cov (  \varphi_k (X_0), X_i - \varphi_k (X_i) )     | \, .
\]
Applying again Lemma \ref{lmapourcov} of the Appendix,   we first write that 
\begin{multline*}
 \sum_{i > \ell_k} |\cov (  \varphi_k (X_0), X_i - \varphi_k (X_i) )     |  \\  \ll \sum_{\ell=0}^{\infty}   \Vert P_{0} ( \varphi_k (X_\ell)) \Vert_2 \sum_{i \geq  [2^{-1} \ell_k] +1} 
i^{-1/2} \Vert \E (g_k ( X_i) | {\mathcal F}_0) - \E (g_k ( X_i) ) \Vert_2 
\, .
\end{multline*}
By condition  \eqref{condondeltaninfty} and since $q>1/2$, 
\[
\sum_{\ell \geq 0} \Vert P_{0} ( \varphi_k (X_\ell)) \Vert_2 \ll \sum_{\ell \geq 0}  ( \ell +1)^{-1/2}  \Vert  \E (  \varphi_k (X_j) | V_0) -  \E (  \varphi_k (X_j))  \Vert_2  \leq C \, . 
\]
So, taking into account \eqref{boundnorm2gk} and the fact that $q>1/2$, 
\begin{multline}  \label{comp3sigma2avechat}
 \sum_{i > \ell_k} |\cov (  \varphi_k (X_0), X_i - \varphi_k (X_i) )     |  \\
 \ll   3^{- k (p-1)/(2p)} \ell_k^{ ( 1-q) /2} {\bf 1}_{q >1} + 3^{- k (p-1) (2q-1)/(2p q)}  {\bf 1}_{q <   1} +  k \, 3^{- k (p-1)/(2p)} {\bf 1}_{q=1}  \, .
\end{multline}
Considering the upper bounds \eqref{comp1sigma2avechat}, \eqref{comp2sigma2avechat} and  \eqref{comp3sigma2avechat}, we then derive 
\begin{multline*}
\Big | \sigma^2 -  (  {\hat c}_{k,0} + 2  \sum_{i \geq 1}   {\hat c}_{k,i} )  \Big | \leq o ( \ell_k 3^{-k (p-2)/p}) 
+ \ell_k^{1/2 - q }   3^{- k (p-2)/(2p)}  +  \big ( \ell_k^{1-  2q } + 3^{- k (p-1)/(2p)} \ell_k^{ ( 1-q) /2} \big )  {\bf 1}_{q> 1}  \\ +  k    3^{- k (p-1)/(2p)} {\bf 1}_{q=1}  
 + \big ( k \,  3^{- k (p-1)/(2p)} + \ell_k^{ 1- 3 q/2  }   3^{- k (p-1)/(2p)} +  3^{- k (p-1) (2q-1)/(2p q)}  \big )  {\bf 1}_{q< 1} 
 \, ,
\end{multline*}
which combined with \eqref{compnukavechat} gives
\begin{multline} \label{compnukavecsigma2}
\big | \nu_k -  \sigma^2 \big |  \ll  o ( \ell_k 3^{-k (p-2)/p}) + \ell_k^{1/2 - q }   3^{- k (p-2)/(2p)}  +  \big ( \ell_k^{1-  2q } + 3^{- k (p-1)/(2p)} \ell_k^{ ( 1-q) /2} \big )  {\bf 1}_{q> 1}
\\  + k    3^{- k (p-1)/(2p)} {\bf 1}_{q=1}   + \big ( k \,  3^{- k (p-1)/(2p)} + \ell_k^{ 1- 3 q/2  }   3^{- k (p-1)/(2p)} +  3^{- k (p-1) (2q-1)/(2p q)}  \big )  {\bf 1}_{q< 1}  \\
+ 3^{k/(2p)} m_k^{-q/2} + \big (  m_k^{-1/2} (\log m_k) {\bf 1}_{q =1}  + m_k^{-q +1/2} {\bf 1}_{q <1} \big )  \, .
\end{multline}
Let us verify that \eqref{cond2v_kbis} holds, namely: 
\[
3^{k(p-2)/(2p)} (  \nu_k  - \sigma^2) =  o((\log k)^{-1/2}) \, .  
\]
The choice of $\ell_k$ implies that $\ell_k 3^{-k (p-2)/p} = 3^{-k(p-2)/(2p)}  (\log k)^{-1/2}$ and $\ell_k^{1/2 - q }  = o((\log k)^{-1/2}) $ (since $q>1/2$). Moreover, 
when $q >1$, we clearly have $ 3^{k(p-2)/(2p)}  \ell_k^{1-2q} =  o((\log k)^{-1/2}) $ and  $3^{k(p-2)/(2p)}  3^{- k (p-1)/(2p)} \ell_k^{ ( 1-q) /2} =  o((\log k)^{-1/2}) $.  It is also clear that 
$3^{k(p-2)/(2p)} k    3^{- k (p-1)/(2p)} =  o((\log k)^{-1/2})$. Next, since $q >(p-1)/2$, 
\[
3^{k(p-2)/(2p)} 3^{- k (p-1) (2q-1)/(2p q)} {\bf 1}_{q< 1}  \leq 3^{k(p-2)/(2p)} 3^{- k (p-1) (p-2)/(2p)} \, ,
\]
proving (since $p>2$) that $ 3^{k(p-2)/(2p)} 3^{- k (p-1) (2q-1)/(2p q)} {\bf 1}_{q< 1}  =  o((\log k)^{-1/2}) $. Also, since $p>2$, 
\[
3^{k(p-2)/(2p)}  \ell_k^{ 1- 3 q/2  }   3^{- k (p-1)/(2p)}  {\bf 1}_{( p-1)/2 < q< 1} \ll  3^{k(p-2)/(2p)}  \ell_k^{ 1/4 }   3^{- k (p-1)/(2p)}  {\bf 1}_{ 2 < p <3}  \, ,
\]
which proves that $3^{k(p-2)/(2p)}  \ell_k^{ 1- 3 q/2  }   3^{- k (p-1)/(2p)}  {\bf 1}_{( p-1)/2 < q< 1} =  o((\log k)^{-1/2})$.
Next, we note that 
\[
 3^{k(p-2)/(2p)}  \frac{3^{k/(2p)} }{m_k^{q/2}} =  3^{k(p-1)/(2p)} 3^{- q k (1-\varepsilon)/p}  = o((\log k)^{-1/2}) \, ,
\]
since $\varepsilon < 1 - \frac{p-1}{ 2q} $.  

Now, if $q =1$ then  $p < 3$ (since $q >(p-1)/2$). Hence since $\varepsilon < 1/2$, we get that $3^{k(p-2)/(2p)} m_k^{-1/2} (\log m_k) {\bf 1}_{q =1}  =  o((\log k)^{-1/2})$. Finally, using again  that $q >(p-1)/2$ and that $\varepsilon < 1/2$, we derive that $3^{k(p-2)/(2p)} m_k^{-q +1/2} {\bf 1}_{q <1}  =  o((\log k)^{-1/2})$.
This ends the proof of \eqref{cond2v_kbis} and then of the theorem. \hfill $\square$

\subsection{Proof of Theorem \ref{KMTavecdeltan1}}  
By Remark \ref{remarkequivcond}, we know that  condition \eqref{condonQdelta-2} is equivalent to \eqref{notaforR}, namely: 
\[
\int_0^1 R^{p-1} (u) Q(u) du < \infty \, , 
\] 
where, for any $u \in [0,1]$, $\gamma^{-1} (u) =\delta^{-1} \circ H(u)$ and  $R(u) = \gamma^{-1} (u)  Q(u)$.

Notice first that, by Proposition 1 in Dedecker-Doukhan \cite{DD03},
\[
\sum_{i \geq 1} |{\rm Cov} (X_0,X_i) | \leq 2 \sum_{i \geq 1} \int_0^{2^{-1} \Vert \E (X_i | V_0) - \E (X_i) \Vert_1 } Q \circ H^{-1} (u) du \leq \int_0^1 R(u) Q(u) du < \infty \, ,
\]
by condition \eqref{notaforR}. Hence the series $\sigma^2 =  {\rm Var} (X_1) + 2 \sum_{i \geq 1} {\rm Cov} (X_1, X_{i+1} )$ converge absolutely and condition 
\eqref{defsigma2asalimit} of Proposition \ref{generalpropBLW} 
holds. 

\medskip

Assume first that $\sigma^2 >0$. To prove the theorem, we shall verify that the other conditions of Proposition \ref{generalpropBLW} are satisfied and with this aim we need to define suitable sequences 
$(m_k)$ and $(M_k)$.  Since we have ${\rm Var}(S_n) /n \rightarrow \sigma^2 >0$, 
it follows that  ${\rm Var}(S_n)  \rightarrow \infty$. Hence $\p (|X_1| >0) >0$ since otherwise we would have $X_1 =0$ a.s. and then $S_n =0$ a.s. for all $n \geq 1$, contradicting the fact that 
${\rm Var}(S_n)  \rightarrow \infty$.  Let $u_1 = (1/2) \p (|X_1| >0)$ (hence $u_1 >0$) and define 
\[
K_0 = \inf \{ k \in  {\mathbb N}\, : \, R(u_1) \leq 3^{k/p} \} \ .
\]
Obviously $K_0 < \infty $ since $u_1>0$ which implies that  $Q(u_1) < \infty  $ and $\gamma^{-1} (u_1) < \infty$. Next,  for any $ k  \geq K_0$, let 
\[
v_k = \inf \{ u \in [0, u_1] \, : \, R(u) \leq 3^{k/p} \} := R^{-1} (3^{k/p}) \, ,  \, M_k = Q(v_k)     \, ,
\]
and $M_k =1$ for $0 \leq k < K_0$. 
Since $u_1 < \p (|X_1| >0) $, it follows that  $Q(u_1) >0$ and therefore since $Q$ is non-increasing and $v_k \leq u_1$, $M_k \geq Q(u_1) > 0$, for $k \geq K_0$. Let now,  for any $ k  \geq K_0$,
\[
m_k = \inf \{n \geq 0 \, : \,  \gamma (n) \leq v_k \}  \, ,
\]
and $m_k = 1$ for any $1 \leq k < K_0$. 
Since $v_k$ is assumed to be strictly less than $1$ (since $v_k \leq u_1 \leq 1/2$), $m_k \geq 1$ (indeed $\gamma (0 )= H^{-1} (\E(  |X_1|) )=1$). In addition, since $R$ is right continuous and non-increasing, $u < R^{-1} (x) \iff R(u) > x$. Hence, $R(R^{-1} (u) ) \leq u $ for all $u \in [0,1]$, implying that 
\[
m_k M_k \leq R(v_k)  \leq 3^{k/p} \, .
\]
Therefore, for any $k \geq K_0$, since $M_k \geq Q(u_1) > 0$, 
\[
m_k \leq (Q(u_1))^{-1} m_k M_k  \leq 3^{k/p} \, ,
\]
which proves that $m_k =o ( 3^{2k/p} k^{-1}) $, as $k \rightarrow \infty$. 

To prove now that the conditions \eqref{condtruncature}, \eqref{condmkdependence}, \eqref{condtoapplySakhanenko} and \eqref{cond2v_k} of Proposition \ref{generalpropBLW} are satisfied, we first notice the following useful facts: 
\beq \label{factonquantile}
Q_{ |\varphi_k (X_j)|} (x) := {\bar Q}_k (u) = Q( x  \vee v_k)  \ \text{ and } \ Q_{ | g_k ( X_j )|} (u) := {\tilde Q}_k (u) =  Q(x) {\bf 1}_{x \leq v_k}\, .
\eeq
Let us start by proving that condition \eqref{condtruncature} holds. By using \eqref{factonquantile}, we get
\[
\sum_{k \geq 1} 3^{k (p-1)/p}  \E (|g_k (X_1)|) = \sum_{k \geq 1} 3^{k (p-1)/p}  \int_0^1 Q_{|g_k (X_1)|} (u) du =  \sum_{k \geq 1} 3^{k (p-1)/p}  \int_0^1 Q (u)  {\bf 1}_{u < v_k} du  \, .
\]
But
\beq \label{evidentpourbornesum}
 \sum_{k \geq 1} 3^{k (p-1)/p}  \int_0^1 Q (u)  {\bf 1}_{u < v_k} du =\int_0^1 Q(u)  \sum_{k \geq 1} 3^{k (p-1)/p}    {\bf 1}_{R(u) > 3^{k/p}} \ll  \int_0^1 R^{p-1} (u) Q(u) du < \infty \, , 
\eeq
by condition \eqref{notaforR} (which is equivalent to condition \eqref{condonQdelta-2}). Hence condition \eqref{condtruncature} is satisfied. Next we note that by Lemma \ref{lmaenlevertildedeltan} of the Appendix, 
\[
\Big \Vert  \max_{1 \leq \ell \leq 3^k - 3^{k-1}} \big |   { W}_{k, \ell} -  {\widetilde W}_{k, \ell} \big | \Big \Vert_1 \leq  \sum_{i=1+3^{k-1}}^{ 3^{k}} \Vert X_{k,i}  - {\tilde X}_{k,i}  \Vert_1 \leq   2 \times  3^k \delta (m_k)  \, .
\]
Therefore, by using \eqref{evidentpourbornesum}, 
\begin{multline} \label{condmkdependenceverify}
\sum_{k \geq 1} 3^{-k/p} \Big \Vert  \max_{1 \leq \ell \leq 3^k - 3^{k-1}} \big |   { W}_{k, \ell} -  {\widetilde W}_{k, \ell} \big | \Big \Vert_1  \leq  2 \sum_{k \geq 1} 3^{k (p-1)/p}  \delta (m_k)  \\
  \leq  2 \sum_{k \geq 1} 3^{k (p-1)/p} H (  \gamma(m_k) ) \leq  2 \sum_{k \geq 1} 3^{k (p-1)/p} \int_0^{v_k } Q(u) du < \infty \, .
 \end{multline}
Hence, condition \eqref{condmkdependence} is satisfied with $\alpha =1$. We prove now that  we can find a real number $ r \in ]2, \infty [$ such that 
 \eqref{condtoapplySakhanenko} holds. With this aim we start by noticing that, for any $r\geq 1$,   by Lemma \ref{lmaenlevertildedeltan} of the Appendix, 
\begin{multline*}
\Big \Vert  \max_{1 \leq \ell \leq  3m_k }  \big |   {\widetilde W}_{k, \ell}   -  {W}_{k, \ell}     \big | \Big  \Vert_r \leq  \sum_{i=3^{k-1} +1}^{3m_k +3^{k-1} }  \Vert X_{k,i}  - {\tilde X}_{k,i}  \Vert_{r}  \\ \leq  (2 M_k)^{(r-1)/r}   \sum_{i=3^{k-1} +1}^{3m_k +3^{k-1} } 
\Vert X_{k,i}  - {\tilde X}_{k,i}  \Vert_1^{1/r}  \ll m_k M_k^{(r-1)/r}  
(\delta (m_k))^{1/r} \, .
\end{multline*}
 Hence, since $m_k M_k \leq 3^{k/p} $, for any $r \geq 1$, 
\[
 \sum_{k \geq \ell_0} \frac { 3^{k} }{  3^{kr /p}m_k}  \E \Big (  \max_{1 \leq \ell \leq  3m_k }  \big |   {\widetilde W}_{k, \ell}   -  {W}_{k, \ell}   \big |^r \Big ) \ll  \sum_{k \geq \ell_0} \frac { (m_k M_k)^{r-1}}{  3^{k (r -p) /p}} \delta (m_k)  
 \ll  \sum_{k \geq \ell_0} 3^{k (p-1)/p} \delta (m_k) \, , 
\]
which is finite by taking into account  \eqref{condmkdependenceverify}.  Hence to prove that condition  \eqref{condtoapplySakhanenko} holds, it suffices to prove that  we can find a real number $r \in ]2, \infty[$ such that
\beq \label{butwithr}
\sum_{k \geq \ell_0} \frac { 3^{k} }{  3^{kr /p}m_k}  \E \Big (  \max_{1 \leq \ell \leq  3m_k }  \big |  W_{k, \ell}    \big |^r \Big ) < \infty \, .
\eeq

To prove \eqref{butwithr}, we  apply the Rosenthal inequality for $\tau$-dependent sequences as given in Corollary 1 in Dedecker-Prieur \cite{DP04}. Let us first recall the definition  of the $\tau$-dependence coefficients:
for any random variable $Y$ with values in ${\mathbb R}^{\ell}$ and any $\sigma$-algebra ${\mathcal F}$, 
\[
\tau \big  ( {\mathcal F} ,  Y  )  = \sup_{f \in {\Lambda_1}({\mathbb R}^{\ell})} \Big  \Vert \sup_{h  \in {\Lambda_1}({\mathbb R})}  \Big \vert  \int h(x) P_{ f(Y) | {\mathcal F}} (dx) -  \int h(x) P_{ f(Y)} (dx) \Big \vert  \Big \Vert_1 \, ,
\]
where, for any integer $\ell \geq 1$, ${\Lambda_1}({\mathbb R}^{\ell})$ is the set of $1$-Lipschitz function from ${\mathbb R}^{\ell}$ to ${\mathbb R}$ with respect to the norm 
$|x-y|_1 \leq \sum_{k=1}^{\ell} |x_i -y_i|$. 
Taking  ${\mathcal F}_p = \sigma ( X_i, i \leq p ) $, the coefficients 
$\tau(i)$ of the sequence $(\varphi_k(X_i))_{i \in {\mathbb Z}}$ are
then defined by:
for any $i \geq 0$, 
\[
\tau(i) = \sup_{u >0} \max_{ 1 \leq \ell \leq u} \frac{1}{\ell} \sup \big \{ \tau \big  ( {\mathcal F}_p ,  ( \varphi_k (X_{j_1}), \ldots,  \varphi_k (X_{j_\ell}))  \big )  \, , \, p+ i \leq j_1 < \ldots < j_{\ell} \big \} \, . 
\]
In the stationary case, Corollary 1 in Dedecker-Prieur \cite{DP04} implies  that,  for any $r >2$, 
\begin{multline*}
 \E \Big (  \max_{1 \leq \ell \leq  3m_k }  \big |  W_{k, \ell}    \big |^r \Big ) \ll \Big ( m_k \int_0^{\Vert \varphi_k (X_1) \Vert_1} (  (\tau/2)^{-1} (u) \wedge m_k) Q_{|\varphi_k (X_1)|} \circ H^{-1}_{|\varphi_k (X_1)|}  (u) du \Big )^{r/2} \\
  + m_k \int_0^{\Vert \varphi_k (X_1) \Vert_1} (   (\tau/2)^{-1} (u) \wedge m_k)^{r-1}Q^{r-1}_{|\varphi_k (X_1)|} \circ H^{-1}_{|\varphi_k (X_1)|}  (u) du \, ,
\end{multline*}
where ${\tau}^{-1}$ is the generalized inverse of the function $\tau$ defined by $\tau (x) = \tau([x])$.


To compare the coefficients $\tau(i)$ with the coefficients $\delta (i)$, we consider $(W_0', ( \varepsilon_j')_{j \geq 1})$  an independent copy of $(W_0, ( \varepsilon_j)_{j \geq 1})$ and define $W_1' = F ( \varepsilon'_1, W_0')$ and $W_m' = F ( \varepsilon_m, W_{m-1}' )$ for any $m \geq 2$. 
Note that for any $j \geq 2$, by using  the relation \eqref{anotherwayforWj} of the Appendix, we have
\beq \label{defXjavecW1}
X_j = h ( \varepsilon_j , W_{j-1}):=h( \varepsilon_j, F_{j-2} (  \varepsilon_{j-1}, \ldots, \varepsilon_2, W_1) ):=g_{j-2}  ( \varepsilon_j,  \varepsilon_{j-1}, \ldots, \varepsilon_2, W_1) \, .
\eeq
Define now, for any $j \geq 2$, 
\beq \label{defXjstaravecW1}
X^{\prime}_j = h ( \varepsilon_j , W'_{j-1}) = h ( \varepsilon_j, F_{j-2} (  \varepsilon_{j-1}, \ldots, \varepsilon_2, W'_1) ):=g_{j-2}  ( \varepsilon_j,  \varepsilon_{j-1}, \ldots, \varepsilon_2, W'_1) \, .
\eeq
Clearly for any $ 2 \leq  j_1 < \ldots < j_{\ell}$, $( \varphi_k (X^{\prime}_{j_1}), \ldots,  \varphi_k (X^{\prime}_{j_\ell})  )$ is distributed as $( \varphi_k (X_{j_1}), \ldots,  \varphi_k (X_{j_\ell})  )$ and is independent of 
$(\varepsilon_0, W_{-1})$. Hence, by stationarity and Lemma 3 in Dedecker-Prieur \cite{DP04}, 
\begin{multline*}
 \tau  \big  ( {\mathcal F}_0 ,  ( \varphi_k (X_{j_1}), \ldots,  \varphi_k (X_{j_\ell})  ) \big )  \\  \leq   \sup_{f \in {\Lambda_1}({\mathbb R}^{\ell})} \Vert  f (\varphi_k (X_{j_1}), \ldots,  \varphi_k (X_{j_\ell})  )   - f ( \varphi_k (X^{\prime}_{j_1}), \ldots,  \varphi_k (X^{\prime}_{j_\ell})  )\Vert_1 \\
  \leq \sum_{v=1}^{\ell} \Vert X_{j_v} - X^{\prime}_{j_v} \Vert_1 = \sum_{v=1}^{\ell}  \iint \E |X_{j_v -1,x} -X_{j_v -1,y} | d\nu(x) d\nu (y) \leq 2 \sum_{v=1}^{\ell}  \delta (j_v ) \, , 
\end{multline*}
where the second inequality comes from the fact that $f \in {\Lambda_1}({\mathbb R}^{\ell})$ and $\varphi_k$ is $1$-Lipschitz. Therefore, since $\delta$ is non-increasing, for any $i \geq 2$,
\[
\tau (i) \leq 2 \delta (i) \, .
\]
Moreover, for any $i \in \{0,1\}$, we obviously get that $\tau (i) \leq 2 \E (|X_1|) = 2 \delta (0)$. It follows that for any $x \geq 0$,
\[
\tau (x):= \tau([x]) \leq 2 \delta (x) \, .
\]
Therefore, since both $\tau$ and $\delta$ are non-increasing, 
\[
(\tau/2)^{-1} (u) \leq \delta^{-1} (u) \, .
\]
In addition,  since $\varphi_k$ is $1$-Lipschitz and such that $\varphi_k(0) =0$, 
\beq \label{relationforH}
H_{|\varphi_k (X_1)|}  ( x) \leq H  ( x) \text{ and then } H^{-1} (u)  \leq H^{-1}_{|\varphi_k (X)|}  (u) \, ,
\eeq
since $H$ is non-decreasing. 
Therefore, using additionally the fact that $u < v \iff Q_{|\varphi_k (X_1)|} (v) < Q_{|\varphi_k (X_1)|} (u)$, we get
\begin{multline*}
 \E \Big (  \max_{1 \leq \ell \leq  3m_k }  \big |  W_{k, \ell}    \big |^r \Big ) \ll \Big ( m_k \int_0^{\Vert X_1\Vert_1} (  \delta^{-1} (u) \wedge m_k) Q_{|\varphi_k (X_1)|} \circ H^{-1} (u) du \Big )^{r/2} \\
  + m_k \int_0^{\Vert X_1 \Vert_1} (   \delta^{-1} (u) \wedge m_k)^{r-1}Q^{r-1}_{|\varphi_k (X_1)|} \circ H^{-1}  (u) du \, ,
\end{multline*}
and then, since $\gamma^{-1} (u) =\delta^{-1} \circ H (u)$, 
\begin{multline*}
 \E \Big (  \max_{1 \leq \ell \leq  3m_k }  \big |  W_{k, \ell}    \big |^r \Big ) \ll \Big ( m_k \int_0^{1} (  \gamma^{-1} (u) \wedge m_k) Q_{|\varphi_k (X_1)|}  (u) Q (u) du \Big )^{r/2} \\
  + m_k \int_0^{ 1} (   \gamma^{-1} (u) \wedge m_k)^{r-1}Q^{r-1}_{|\varphi_k (X_1)|} (u)  Q (u) du \, .
\end{multline*}
Recall now that $m_k =  \gamma^{-1} (v_k)$, therefore since $\gamma^{-1}$ is non-increasing,  
\[ \gamma^{-1} (u) \wedge m_k = m_k {\bf 1}_{ u < v_k } +  \gamma^{-1} (u) {\bf 1}_{ u \geq  v_k } \, .
\]
Using also the fact that $Q_{ |\varphi_k (X_1)|} (x) = Q( x  \vee v_k) $,  we get 
\begin{multline*}
 \E \Big (  \max_{1 \leq \ell \leq  3m_k }  \big |  W_{k, \ell}    \big |^r \Big )  \ll   \Big ( m_k Q^{1-2/r}(v_k)   \int_0^{1}  \gamma^{-1}  (u)  Q^{1+2/r} (u) du \Big )^{r/2}  \\
  + m_k (m_k Q (v_k))^{r-1} \int_0^{ 1}  Q (u) {\bf 1}_{ u < v_k } du   +  m_k \int_0^{ 1} R^{r-1} (u)  Q (u) {\bf 1}_{ u \geq  v_k } du  \\ := I_k^{(1)} + I_k^{(2)} + I_k^{(3)}\, .
\end{multline*}
Using the fact that  $m_k Q(v_k)  \leq 3^{k/p}$ and \eqref{evidentpourbornesum}, we get that, for any $r >2$, 
\[
 \sum_{k \geq  \ell_0}   3^{k(p-r)/p} m_k^{-1} I_k^{(2)}   \leq  \sum_{k \geq  \ell_0}   3^{k(p-1)/p}  \int_0^{ 1}  Q (u) {\bf 1}_{ u < v_k } du   < \infty \, .
\]
On another hand,  for any $r >p$, 
\begin{multline*}
 \sum_{k \geq  \ell_0}   3^{k(p-r)/p} m_k^{-1} I_k^{(3)}   \leq  \sum_{k \geq  \ell_0}   3^{k(p-r)/p}  \int_0^{ 1} R^{r-1} (u)  Q (u) {\bf 1}_{ u \geq  v_k } du \\  =  \int_0^{ 1}  \sum_{k \geq  \ell_0}   3^{k(p-r)/p}   
 R^{r-1} (u)   {\bf 1}_{ 3^{k/p } \geq  R(u)  }  Q (u) du 
 \ll  \int_0^1 R^{p-1} Q(u) du < \infty \, ,
\end{multline*}
by condition \eqref{notaforR} (which is equivalent to condition \eqref{condonQdelta-2}). Finally using again that $m_k Q(v_k)  \leq 3^{k/p}$, we derive that, for any $r >2 (p-1)$, 
\begin{multline*}
 \sum_{k \geq  \ell_0}   3^{k(p-r)/p} m_k^{-1} I_k^{(1)}  \leq  \sum_{k \geq  \ell_0}   3^{k(p-r)/p}  \Big ( \big ( m_k Q(v_k) \big )^{1-2/r} \int_0^1 \gamma^{-1} (u)  Q^{1+2/r}( u )  du \Big )^{r/2} \\
 \leq  \sum_{k \geq  \ell_0}   3^{k(p-r)/p}  3^{k (r-2)/(2p)} \Big (  \int_0^1 \gamma^{-1} (u)  Q^{1+2/r}( u )  du \Big )^{r/2}  \, ,
\end{multline*}
since condition \eqref{notaforR} obviously implies that $\int_0^1 \gamma^{-1} (u)  Q^{1+2/r}( u )  du < \infty$. 
So, overall, \eqref{butwithr} holds provided  we select $r > 2 (p-1)$. 

\smallskip

To end the proof it remains to show that condition \eqref{cond2v_k} holds.  With this aim, we start by recalling  the equation \eqref{nuk1}, namely: 
\[
\nu_k = {\tilde c}_{k,0}  + 2 \sum_{\ell =1}^{m_k }  {\tilde c}_{k,\ell} \, , 
\]
where, for $i \geq 0$, 
\[
{\tilde c}_{k,i}= {\rm cov} ({\tilde X}_{k, m_k+1} , {\tilde X}_{k,i+m_k+1} ) \, \text{ and } {\hat c}_{k,i} = {\rm cov} ({ X}_{k,0} , { X}_{k,i} ) \, .
\]
But, by using Lemma \ref{lmaenlevertildedeltan} of the Appendix, we have, for any $i \geq 0$, 
\begin{multline*}
\big | {\tilde c}_{k,i} - {\hat c}_{k,i} \big |  =   \big |   {\rm cov} ({\tilde X}_{k,m_k+1} -{ X}_{k,m_k+1}, {\tilde X}_{k,i+m_k+1} )  +  {\rm cov} ({ X}_{k,m_k+1}, {\tilde X}_{k,i+m_k+1} -{X}_{k,i+m_k+1} )  \big |    \\ 
\leq   Q(v_k)  \Vert {\tilde X}_{k,m_k+1} -{ X}_{k,m_k+1} \Vert_1 +  Q(v_k)  \Vert  {\tilde X}_{k,i+m_k+1} - X_{k,i+m_k+1}   \Vert_1  \\  \leq  4 \times Q(v_k)  \delta (m_k)  \leq  4 \times Q(v_k) \int_0^{\gamma (m_k)} Q(u) du  
\leq  4 \times Q(v_k) \int_0^{v_k } Q(u) du \, .
\end{multline*}
Hence, since $m_k Q(v_k) \leq 3^{k/p}$, 
\begin{multline} \label{nuk2deltan}
\Big  | \nu_k -  {\hat c}_{k,0}  - 2 \sum_{\ell =1}^{m_k }  {\hat c}_{k,\ell}   \Big | \leq   8 \times m_k Q(v_k) \int_0^{1 } Q(u) {\bf 1}_{u < v_k}du \leq 8 \times 3^{k/p}  \int_0^{1 } Q(u) 
 {\bf 1}_{R(u) >3^{k/p}}du \\
 \leq  8 \times 3^{k/p}  3^{ - k (p-1)/p}   \int_0^{1 } R^{p-1}Q(u) du \ll  3^{ - k (p-2)/p} \, ,
\end{multline}
by condition \eqref{notaforR} (which is equivalent to condition \eqref{condonQdelta-2}). Taking into account \eqref{nuk2deltan} together with the fact that $\sigma^2 = \sum_{k \in {\mathbb Z}} {\rm cov} (X_0, X_k) $, we get 
\beq \label{diffnuksigma2inter}
|\nu_k  - \sigma^2| \ll 2   \sum_{i=0}^{m_k} |  {\hat c}_{k,i } -  {\rm cov} (X_0, X_i) | +  2 \sum_{i \geq m_k +1} | {\rm cov} (X_0, X_i) | +  3^{ - k (p-2)/p} \, .
\eeq
Next, by using Proposition 1 in Dedecker-Doukhan \cite{DD03}, we derive
\[
\sum_{i\geq m_k} | {\rm cov} (X_0, X_i) | \leq  2  \sum_{i\geq m_k}  \int_0^{\delta(i)} Q \circ H^{-1} (u) du = 2  \sum_{i\geq m_k}  \int_0^{\gamma(i)} Q^2 (u) du  \, .
\]
But, since $m_k = \gamma^{-1}(v_k)$, note that 
\begin{multline} \label{usefulfact}
\int_0^1 R(u) Q (u) 1_{u < v_k}du = \int_0^1 \gamma^{-1} (u) Q^2 (u)
{\mathbf 1}_{u < v_k} du = \sum_{i  \geq 0} \int_0^1  Q^2 (u) {\mathbf 1}_{u < v_k}  {\mathbf 1}_{u <
\gamma (i) } du \\
= m_k \int_0^1  Q^2 (u) 1_{u < v_k}  du + \sum_{i \geq m_k} \int_0^1  Q^2
(u)  {\mathbf 1}_{u < \gamma (i) } du \, .
\end{multline}
Hence
\begin{multline}  \label{boundrestesommecov}
\sum_{i\geq m_k} | {\rm cov} (X_0, X_i) |  \leq 2 \int_0^1 R(u) Q (u) 1_{u < v_k}du = 2 \int_0^1 R(u) Q (u) 1_{R(u) > 3^{k/p} }du \\
\leq 2 \times  3^{ - k (p-2) /p} \int_0^1 R^{p-1}(u) Q (u) 1_{R(u) > 3^{k/p} }du \ll 3^{ - k (p-2) /p} \, ,
  \end{multline}
by condition \eqref{notaforR}. On another hand,  by using inequality (1.11a) in \cite{Ri00} and \eqref{factonquantile}, we derive that, for any $i \geq 0$,
 \begin{multline*}
  | {\hat c}_{k,i} - {\rm cov} (X_0, X_i) | =  \big |   {\rm cov} (X_{k,0} -{ X}_{0}, X_{k,i} )  +  {\rm cov} (X_{0}, X_{k,i} -{X}_{i} )  \big |    \\ 
  \leq  2 \int_0^{1/2} Q_{|\varphi_k (X)|} (u) Q_{|g_k (X)|} (u) du  + 2 \int_0^{1/2} Q_{|X_1|} (u) Q_{|g_k (X)|} (u) du \leq  4   \int_0^1  Q^2 (u)  {\bf 1}_{ u < v_k}du \, .
  \end{multline*}
Hence, by taking into account \eqref{usefulfact}, 
  \[
  \sum_{i=0}^{m_k} |  {\hat c}_{k,i } -  {\rm cov} (X_0, X_i) | \leq 4 (m_k +1)  \int_0^1  Q^2 (u)  {\bf 1}_{ u < v_k}du \leq 8 \int_0^1 R(u) Q (u) 1_{u < v_k}du  \, .
  \]
  So, by the computations in \eqref{boundrestesommecov}, 
  \beq  \label{bounddiffbetweencov}
   \sum_{i=0}^{m_k} |  {\hat c}_{k,i } -  {\rm cov} (X_0, X_i) | \ll 3^{ - k (p-2) /p} \, .
  \eeq
 Hence, starting from \eqref{diffnuksigma2inter} and taking into account \eqref{boundrestesommecov} and \eqref{bounddiffbetweencov}, it follows that 
 \[
 | \nu_k  - \sigma^2 | \ll    3^{ - k (p-2)/p}   \, ,
 \]
implying, since $p>2$,  that 
\[
3^k  (  \nu_k - \sigma^2)^2  \ll     3^{ k (4-p)/p}  = o(3^{2k/p} (\log k)^{-1})  \, , \, \mbox{ as $k \rightarrow \infty$}\, .  
\]
This proves that  \eqref{cond2v_kbis} holds and then  that \eqref{cond2v_k} is satisfied since $\sigma^2>0$.  The proof is complete for the case $\sigma^2>0$.  

\medskip

Assume now that $\sigma^2 =0$.  Let $M$ be a positive real number. According to inequality (5.42) in Merlev\`ede-Rio \cite{MR12}, for any positive integer $n$,  any real number $\lambda$, and 
any positive integer $q \leq n$ and such that $qM \leq \lambda$, we have 
\begin{multline} \label{decproba}
 \p \Big ( \max_{1 \leq k \leq n} |S_k -k \E (X_1)| \geq 5 \lambda \Big ) \leq  \frac{{\rm Var}(S_n) }{ \lambda^2} + 2 \lambda^{-1}\sum_{i=1}^n \E | g_M (X_i)| \\ +
\lambda^{-1}\sum_{i=1}^n  \| \E( \varphi_M (X_i) | V_{i-q}) - \E ( \varphi_M (X_i)) \|_1 \, .
  \end{multline}
Choose now $u = R^{-1}(\lambda)$, $q = \gamma^{-1} (u) \wedge n $ and $M = Q (u)$. Since $R$ is right
continuous, we have $R(u)\leq \lambda$, hence $qM \leq R(u)  \leq \lambda$. Note also that 
  \begin{equation} \label{dec13FN}
  \sum_{k=1}^n  \E | g_M (X_k)| 
  \leq   n \int_0^u  Q (x) dx
  \leq  n \int_0^1 Q(x) {\bf 1}_{R(x)> \lambda} dx \, .
  \end{equation}
In addition, 
\[
 \| \E( \varphi_M (X_i) | V_{i-q}) - \E ( \varphi_M (X_i)) \|_1 \leq 2 \delta (q) = 2 \int_0^{ \gamma (q)} Q(x) dx  \, .
\]
Since $\gamma(q)  \leq u$, it follows that 
\[
\sum_{i=1}^n  \| \E( \varphi_M (X_i) | V_{i-q}) - \E ( \varphi_M (X_i)) \|_1  \leq  2 n  \int_0^1 Q(x) {\bf 1}_{R(x)> \lambda}  dx \, .
\]
Starting from \eqref{decproba} and taking into account the considerations above, we get that, for any $\lambda >0$,
\[
 \p \Big ( \max_{1 \leq k \leq n} |S_k -k \E (X_1)| \geq 5 \lambda
 \Big ) \leq  \frac{{\rm Var}(S_n) }{ \lambda^2} +  \frac{ 4 n}{\lambda} \int_0^1 Q(x) {\bf 1}_{R(x)> \lambda} dx \, .
\]
Hence, for any $\varepsilon >0$, selecting $\lambda = \varepsilon n^{1/p}$, we derive
\begin{multline*} 
\sum_{n \geq 1} n^{-1} \p \Big ( \max_{1 \leq k \leq n} |S_k -k \E (X_1)| \geq 5 \varepsilon n^{1/p}
\Big ) \\  \leq  \varepsilon^{-2} \sum_{n \geq 1}   \frac{{\rm Var}(S_n) }{  n^{1 + 2/p} } +  4 \varepsilon^{-1}  \sum_{n \geq 1} n^{-1/p}
\int_0^1 Q(x) {\bf 1}_{R(x)>  \varepsilon n^{1/p}} dx   \, .
  \end{multline*}
  The second series in the right-hand side is finite under condition condition \eqref{notaforR} (which is equivalent to condition \eqref{condonQdelta-2}). 
Hence, if we can prove that 
\beq \label{sufficientcondvarsigmenum}
\sum_{n \geq 1}   \frac{{\rm Var}(S_n) }{  n^{1 + 2/p} }  < \infty\, , 
\eeq
then we will get that,  for any $\varepsilon >0$, 
\[
\sum_{n \geq 1} n^{-1} \p \Big( \max_{1 \leq k \leq n} |S_k -k \E (X_1)| \geq  \varepsilon n^{1/p}\Big ) < \infty 
\]
which will imply  $
S_n - n \E (X_1) = o (n^{1/p})$ a.s. and therefore the proof of the theorem will be complete. In the case where $p \geq 3$, \eqref{sufficientcondvarsigmenum} is almost immediate. To see this, we first note that  condition \eqref{notaforR} implies 
$\sum_{i \geq 1} i  |{\rm Cov} (X_0,X_i) | < \infty$. Indeed, by Proposition 1 in Dedecker-Doukhan \cite{DD03},
\[
\sum_{i \geq 1} i |{\rm Cov} (X_0,X_i) | \leq 2 \sum_{i \geq 1} i  \int_0^{2^{-1} \Vert \E (X_i | V_0) - \E (X_i) \Vert_1 } Q \circ H^{-1} (u) du \leq \int_0^1 R^2(u) Q(u) du \, , 
\]
which is finite under condition \eqref{notaforR}. Therefore, by Lemma 1 in Bradley \cite{Br97}, ${\rm Var}(S_n)$ is bounded which obviously entails \eqref{sufficientcondvarsigmenum}. To handle the case where $p \in ]2,3[$, we first note that, by inequality (4.84) in \cite{MP06}, 
\[
\Vert \E ( X_k | V_0) - \E (X_k)\Vert^2_2 \leq  \E \big | X_k ( \E ( X_k | V_0) - \E (X_k) )\big | \leq 3 \int_0^{\Vert  \E ( X_k | V_0) - \E (X_k) \Vert_1 }  Q \circ H^{-1} (u) du \, . 
\]
But, $\Vert  \E ( X_k | V_0) - \E (X_k) \Vert_1 \leq 2 \delta (k) $. Hence 
\[
\Vert \E ( X_k | V_0) - \E (X_k)\Vert^2_2 \leq 6 \int_0^{ \delta(k) }  Q \circ H^{-1} (u) du \, .
\]
Hence condition \eqref{condonQdelta-2} entails
\[
\Vert \E ( X_k | V_0) - \E (X_k)\Vert_2 \ll k^{-(p-1)/2} \, , 
\]
which implies (since $p>2$) that 
\beq \label{toapplythemartapprox}
\sum_{k \geq 1} k^{-1/2} \Vert \E ( X_k | V_0) - \E (X_k)\Vert_2 < \infty   \, \mbox{ and therefore } \, \sum_{k \geq 1}  \frac{ \Vert \E ( S_k | V_0) -  k \E (X_1)\Vert_2}{k^{3/2}} < \infty \, .
\eeq
We use now the same arguments as developed at the beginning of the proof of Theorem  \ref{KMTavecdeltaninfty}. The fact that the series in \eqref{toapplythemartapprox} converge 
implies that there exists a stationary sequence $(D_k)_{k \in {\mathbb Z}}$ of martingale differences in ${\mathbb L}^2$, adapted to 
$({\mathcal F}_k)_{k \in {\mathbb Z}}$ and such that 
$$
\lim_{n \rightarrow \infty} n^{-1/2} \Big \Vert S_n - n \E (X_1) - \sum_{k=1}^n D_k \Big \Vert_2 = 0 \, .
$$
 Together with the fact that 
 $\lim_{n \rightarrow \infty} n^{-1}{\rm Var} (S_n) = \sigma^2 =0$, it follows that $D_k =0$ a.s, for any $k$.  Hence, using  the upper bound (4) in Cuny-Merlev\`ede \cite{CM14} (see also Proposition 1 in \cite{MPP12}), it follows that, for any $p \in ]2,3[$, 
\begin{multline*}
{\rm Var} (S_n) \ll n \Big (  \sum_{k \geq n}  \frac{ \Vert \E ( S_k | V_0) -  k \E (X_1)\Vert_2}{k^{3/2}} \Big )^2 \\
\ll n \Big (   \frac{1}{\sqrt n} \sum_{k=1}^n \Vert \E ( X_k | V_0) - \E (X_k)\Vert_2 + \sum_{k \geq n}  \frac{ \Vert \E ( X_k | V_0) -   \E (X_1)\Vert_2}{k^{1/2}} \Big )^2 \ll n^{3-p} \, .
\end{multline*}
Therefore, for any $p \in ]2,3[$, 
\[
\sum_{n \geq 1}   \frac{{\rm Var}(S_n) }{  n^{1 + 2/p} } \ll \sum_{n \geq 1}   \frac{ 1}{  n^{ p+ 2/p -2} } \, , 
\]
which is finite since $ p+ 2/p -3 = p^{-1} (p-1) (p-2) >0$. This ends the proof of the theorem. 
\hfill $\square$

\section{Appendix}

\subsection{Some technical results}

In this section, we collect some technical results that are useful for the proofs of Theorems \ref{KMTavecdeltaninfty} and \ref{KMTavecdeltan1}.

\begin{Lemma} \label{lmapourcov} 
Let $(Y_k)_{k \in {\mathbb Z}}$ be a stationary sequence of real-valued random variables adapted to an increasing and stationary filtration $({\mathcal F}_k)_{k \in {\mathbb Z}}$. Let $f$ and $g$ be two functions in 
${\mathbb L}^2 ({\mathbb R}, P_{Y_0})$ such that $\E (f (Y_0 ) | {\mathcal F}_{- \infty}) = \E (f (Y_0 ) ) $ a.s. and $\E (g (Y_0 ) | {\mathcal F}_{- \infty}) = \E (g (Y_0 ) ) $ a.s. Then, for any positive integer $L$, 
\[
\sum_{i \geq L} | {\rm cov} (f(Y_0), g(Y_i))|  \leq  3 \sqrt{2} \sum_{\ell=0}^{\infty} \Vert P_{0} (f(Y_\ell)) \Vert_2  \Big (  \sum_{k  \geq [(L+\ell) /2] +1} k^{-1/2}  \Vert  \E ( g(Y_k) | {\mathcal F}_0) -  \E ( g(Y_k)  ) \Vert_2  \Big )  
\]
and
\[
\sum_{i =0}^L 
\Vert P_{0} ( g(Y_i)) \Vert_2  \leq   \sqrt{2}  \sum_{k =0}^L (k+1)^{-1/2}  \Vert  \E ( g(Y_k) | {\mathcal F}_0) -  \E ( g(Y_k)  ) \Vert_2  \, , 
\]
where $P_{j} (\cdot) = \E ( \cdot | {\mathcal F}_{j} ) -  \E ( \cdot | {\mathcal F}_{j-1} ) $.
\end{Lemma}
\noindent {\bf Proof.} 
Since $\E (f (Y_0 ) | {\mathcal F}_{- \infty}) = \E (f (Y_0 ) ) $ a.s. and $\E (g (Y_0 ) | {\mathcal F}_{- \infty}) = \E (g (Y_0 ) ) $ a.s., we first write
\[
f(Y_0) - \E ( f(Y_0)) = \sum_{\ell=0}^{\infty}  P_{-\ell} (f(Y_0)) \,  \mbox { and } \,  g(Y_i) - \E ( g(Y_i)) = \sum_{\ell=-i}^{\infty}  P_{-\ell} (g(Y_i))   \mbox { a.s. }
\]
Hence, by orthogonality, for any $i \geq 0$, 
\[
 |  {\rm cov} ( f(Y_0), g(Y_i)) | \leq \sum_{\ell=0}^{\infty}| \E (P_{-\ell} ( f(Y_0)) P_{- \ell} (g(Y_i)) )| \, , 
\]
and then, by Cauchy-Schwarz's inequality and stationarity, 
\beq \label{majcovavecP0}
\sum_{i \geq L} | {\rm cov} (f(Y_0), g(Y_i))| \leq  \sum_{\ell=0}^{\infty} \Vert P_{0} ( f( Y_\ell) ) \Vert_2 \sum_{i \geq  L+\ell} 
\Vert P_{0} ( g(Y_i)) \Vert_2 \, .
\eeq 
But, for any $m \geq 1$,  by Cauchy-Schwarz's inequality, 
\begin{multline} \label{ine1lmalmapourcov}
\sum_{i \geq m} 
\Vert P_{0} ( g(Y_i)) \Vert_2 = \sum_{i \geq  m} i^{-1}
\Vert P_{0} ( g(Y_i)) \Vert_2 \sum_{k =1}^i = \sum_{k  \geq 1} \sum_{i \geq  \max ( m , k) } i^{-1}
\Vert P_{0} ( g(Y_i)) \Vert_2  \\ = \sum_{k = 1}^{m} \sum_{i \geq m  } i^{-1}
\Vert P_{0} ( g(Y_i)) \Vert_2 +  \sum_{k  >m} \sum_{i \geq k  } i^{-1}
\Vert P_{0} ( g(Y_i)) \Vert_2 \\
\leq  \sqrt{ 2 m}   \Big (  \sum_{i \geq m  } 
\Vert P_{0} ( g(Y_i)) \Vert^2_2  \Big )^{1/2} +  \sqrt{ 2 } \sum_{k  >m} k^{-1/2}   \Big (  \sum_{i \geq k   } 
\Vert P_{0} ( g(Y_i)) \Vert^2_2  \Big )^{1/2}  \, , \end{multline}
giving
\[
\sum_{i \geq m} 
\Vert P_{0} ( g(Y_i)) \Vert_2 \leq  \sqrt{2 m}  
\Vert  \E ( g(Y_m) | {\mathcal F}_0) -  \E ( g(Y_m)  ) \Vert_2  + \sqrt{ 2 }  \sum_{k  >m} k^{-1/2}  \Vert  \E ( g(Y_k) | {\mathcal F}_0) -  \E ( g(Y_k)  ) \Vert_2 \, .
\]
Since $(\Vert  \E ( g(Y_k) | {\mathcal F}_0) -  \E ( g(Y_k)  ) \Vert_2)_{k \geq 0}$  is non-increasing, we get that for any $m \geq 1$, 
\[
\sum_{i \geq m} 
\Vert P_{0} ( g(Y_i)) \Vert_2 \leq 3 \sqrt{2}  \sum_{k  \geq [m /2] +1} k^{-1/2}  \Vert  \E ( g(Y_k) | {\mathcal F}_0) -  \E ( g(Y_k)  ) \Vert_2 \, , 
\]
which combined with \eqref{majcovavecP0} gives the first inequality of the lemma. To prove the second one, it suffices to write that  $\sum_{i =0}^L
\Vert P_{0} ( g(Y_i)) \Vert_2 = \sum_{i =0}^L (i+1)^{-1}
\Vert P_{0} ( g(Y_i)) \Vert_2 \big (  \sum_{k =1}^{i+1} 1\big) $ and to use Cauchy-Schwarz's inequality as in \eqref{ine1lmalmapourcov}. 
 \hfill $\square$
 
 \medskip
 
The following proposition is   a  non stationary version of the Peligrad-Utev-Wu \cite{PUW07} inequality. As in  \cite{PUW07}, the proof can be done by induction (a complete proof appears in Section 3.2.1 of \cite{MPU17}).
\begin{Proposition}
\label{LpmainNS} Let $\alpha \geq 2$ and  $(Z_k)_{k \in {\mathbb Z}}$ be a sequence of real-valued random variables in $\LL^{\alpha} $ and adapted to a non-decreasing filtration 
$(\F_k)_{k \in {\mathbb Z}}$. Then, for any $n \geq 1$, 
\begin{multline*}
\Big \Vert  \max_{1 \leq k \leq n} \Big \vert \sum_{i=1}^k Z_i \Big \vert \Big  \Vert_\alpha \leq   ( 2 c_\alpha  + 1) \Big (  \sum_{j=1}^{n} \Vert Z_j \Vert_\alpha^2 \Big )^{1/2}  \\
+  {\sqrt 2}  ( 2 c_\alpha  + 1)
 \sum_{\ell=0}^{r-1}  \Big ( \sum_{m=1}^{2^{r-\ell}} \Vert {\mathbb E} ( S_{m 2^\ell } -S_{(m-1) 2^\ell }  | {\mathcal F}_{(m-2) 2^\ell +1} ) \Vert_\alpha^2 \Big )^{1/2} \, ,
\end{multline*}
where  $S_k = \sum_{i=1}^k Z_i$, $ c_\alpha = \frac{\alpha}{(\alpha-1)^{1/2}}$ if $\alpha >2$, $c_2 = 1$ and $r$ is the unique positive integer such that $2^{r-1} \leq n < 2^r$.  
\end{Proposition}

\begin{Lemma} \label{lmaenlevertildedeltan} 
For any $q \in [1,p )$,   for any $k \geq 1$ and any $j \geq m_k+1$,
\[
\Vert X_{k,j}  - {\tilde X}_{k,j}  \Vert^q_q \leq  \iint  \E (|  X_{m_k+1,x} - X_{m_k+1,y} |^q  ) \nu(dx) \nu (dy) \, ,
\] 
where $X_{k,j}$ and $ {\tilde X}_{k,j} $ are defined in \eqref{defXkj} and  \eqref{defXkjtilde} respectively.
\end{Lemma}
\noindent {\bf Proof.} Let $(W_0', ( \varepsilon_j')_{j \geq 1})$ be an independent copy of $(W_0, ( \varepsilon_j)_{j \geq 1})$ and define $W_j' = F(\varepsilon'_j, W'_{j-1})$, $j \geq 1$. 
For $\ell \geq 1$, let $F_\ell$ be the function from $G^{\ell} \times X$ to $X$ defined in an iterative way as follows  
\[
F_1 = F \text{ and } F_{\ell} (x_1, x_2, \ldots, x_{\ell+1}) = F_{\ell-1} (x_1, x_2, \ldots, x_{\ell-1}, F ( x_{\ell}, x_{\ell+1})) \, , \, \ell \geq 2 \, .
\]
Note that for any integer $\ell $ such that $1 \leq \ell \leq j-1$, 
\beq \label{anotherwayforWj}
W_{j-1} = F_{\ell} (  \varepsilon_{j-1}, \varepsilon_{j-2}, \ldots, \varepsilon_{j - \ell} , W_{j - \ell -1}) \, .
\eeq
Hence, for any $j \geq m_k+1$,  
\begin{multline*}
\E \big ( \varphi_k (X_j)  | \varepsilon_j, \varepsilon_{j-1}, \ldots, \varepsilon_{j-m_k}\big ) \\ = \E \big ( \varphi_k ( h ( \varepsilon_j,  F_{m_k} ( \varepsilon_{j-1}, \ldots, \varepsilon_{j-m_k}, W_{j-m_k-1}) ) )  |  \varepsilon_j, \varepsilon_{j-1}, \ldots, \varepsilon_{j-m_k} \big ) \\
 = \E \big ( \varphi_k ( h ( \varepsilon_j,  F_{m_k} (  \varepsilon_{j-1}, \ldots, \varepsilon_{j-m_k}, W'_{j-m_k-1}) ) )  |  \varepsilon_j, \varepsilon_{j-1}, \ldots, \varepsilon_{1},  W_0 \big ) \, .
\end{multline*}
On another hand, for any $j \geq 1$, 
\begin{multline*}
\varphi_k (X_j)   = \E \big ( \varphi_k (X_j)  |  \varepsilon_j, \varepsilon_{j-1}, \ldots, \varepsilon_{1},  W_0 \big ) \\
 =\E \big ( \varphi_k ( h ( \varepsilon_j,  F_{m_k} (  \varepsilon_{j-1}, \ldots, \varepsilon_{j-m_k}, W_{j-m_k-1}) ) )  |  \varepsilon_j, \varepsilon_{j-1}, \ldots, \varepsilon_{1},  W_0 \big ) \, .
\end{multline*}
Hence,   for any $j \geq m_k+1$,
\begin{multline*}
\Vert X_{k,j}  - {\tilde X}_{k,j}  \Vert_q \\ \leq  \Vert \varphi_k ( h ( \varepsilon_j,  F_{m_k} (  \varepsilon_{j-1}, \ldots, \varepsilon_{j-m_k}, W_{j-m_k-1}) ) )  - \varphi_k ( h ( \varepsilon_j,  F_{m_k} ( \varepsilon_j, \varepsilon_{j-1}, \ldots, \varepsilon_{j-m_k}, W'_{j-m_k-1}) ) )  \Vert_q \\
 \leq  \Vert  h ( \varepsilon_j,  F_{m_k} ( \varepsilon_{j-1}, \ldots, \varepsilon_{j-m_k}, W_{j-m_k-1}) )  - h ( \varepsilon_j,  F_{m_k} ( \varepsilon_{j-1}, \ldots, \varepsilon_{j-m_k}, W'_{j-m_k-1}) )  \Vert_q \, , 
\end{multline*}
where the second inequality comes from the fact that $\varphi_k$ is $1$-Lipschitz. By stationarity, it follows that 
\[
\Vert X_{k,j}  - {\tilde X}_{k,j}  \Vert_q  \leq  
\Vert  h ( \varepsilon_{m_k+1},  F_{m_k} ( \varepsilon_{m_k}, \ldots, \varepsilon_{1}, W_{0}) )  - h ( \varepsilon_{m_k+1},  F_{m_k} ( \varepsilon_{m_k}, \ldots, \varepsilon_{1}, W'_{0}) ) \Vert_q  \, .
\]
Hence, if we define  $(X^*_n)_{n \geq 1} $  by 
\[
X^*_n = h ( \varepsilon_{n}, W^*_{n-1} ) \, \text{ where } W^*_{n} = F ( \varepsilon_{n}, W^*_{n-1} ) \, , 
\]
with  $W^*_{0} $ independent of $(W_0, (\varepsilon_k)_{k \geq 1})$ and such that  $ W^*_{0}=^{\mathcal L}W_{0}$, we get that for any $j \geq m_k+1$,
\beq \label{relcondavecXstar}
\Vert X_{k,j}  - {\tilde X}_{k,j}  \Vert_q \leq  
\Vert X_{m_k+1} - X^{*}_{m_k+1}\Vert_q \, .
\eeq
But,
\begin{align*}
\Vert X_{m_k+1} - X^{*}_{m_k+1}\Vert^q_q & = \iint  \E (|  X_{m_k+1} - X^{*}_{m_k+1}  |^q  | W_0=x, W^*_0=y) \nu (dx) \nu (dy) \\
&  =  \iint  \E (|  X_{m_k+1,x} - X_{m_k+1,y} |^q  ) \nu (dx) \nu(dy) \, , 
\end{align*}
which combined with \eqref{relcondavecXstar} gives the lemma.  \hfill $\square$

\subsection{Proof of Lemma \ref{lmacontinuous}} \label{sectionprooflmacontinuous}

The first inequality in \eqref{aimMC} comes from the coupling inequality \eqref{couplinginebeta} and the fact that $\liminf_{n \rightarrow \infty} n^a \beta (n) >0$ (see Theorem 9.4 in Rio \cite{Ri00}).  We prove now the  second inequality in \eqref{aimMC}. 

Let $W_{n,x}$ be the chain starting at $x$. Note first that for any any $x,y \in [0,1]$, 
\begin{multline*}
\p_{x,y} (T^* >n) = \p_{x,y} \big (T^* >n, \{W_{n,x} =x  \cup W_{n,y} =y \} \big ) +  \p_{x,y} \big (T^* >n, \{W_{n,x} \neq x  , W_{n,y} \neq y \} \big )  \\
\leq (1-x)^n + (1-y)^n  + \p_{x,y} \big (T^* >n, \{W_{n,x} \neq x  , W_{n,y} \neq y \} \big ) \, .
\end{multline*}
But
\begin{multline*}
 \p_{x,y} \big (T^* >n, \{W_{n,x} \neq x  , W_{n,y} \neq y \} \big ) \\
  = \sum_{i=1}^n  \sum_{j=1, j \neq i }^n \p_{x,y} \big (T^* >n,   W_{n,x} =  F_{\pi}^{-1} (V_i)  , W_{n,y} =  F_{\pi}^{-1} (V_j)  \big )  \, .
\end{multline*} 
For $j >i$, define ${\mathcal W}_{i,j} :=  \bigcap_{k=i}^{j} \{W_{k,x} \neq W_{k,y} \} $, ${\mathcal E}_{i,j} (x):=  \bigcap_{k=i}^{j} \{W_{k,x} = F_{\pi}^{-1} (V_i ) \} $, and  note that 
\begin{multline*}
 \p_{x,y} \big (T^* >n,   W_{n,x} =  F_{\pi}^{-1} (V_i)  , W_{n,y} =  F_{\pi}^{-1} (V_j)  \big ) 
\\ = (a+1)^2\iint_{[0,1]^2}  \p_{x,y} \Big ( {\mathcal W}_{1,i-1} ,  {\mathcal E}_{i,n} (x) 
 ,  {\mathcal W}_{i,j-1} ,  {\mathcal E}_{j,n} (y) \big  | F_{\pi}^{-1} (V_i)=u  , F_{\pi}^{-1} (V_j) = v \Big )  u^a v^a du dv  \\
\leq (a+1)^2\iint_{[0,1]^2}  \p_{x,y} \Big ( {\mathcal W}_{1,i}, \, \{U_i<W_{i-1,x} \},\, \bigcap_{k=i +1}^{j} \{  U_k \geq u  \} ,  \bigcap_{k=j +1}^{n} \{  U_k \geq u \vee v   \}  \Big )  u^a v^a du dv \\
\leq (a+1)^2\iint_{[0,1]^2}  \p_{x,y}  ( T^* >i, \, W_{i,x}=F_\pi^{-1}
(V_i)   ) \p \Big  (  \bigcap_{k=i +1}^{j} \{  U_k \geq u  \} ,  \bigcap_{k=j +1}^{n} \{  U_k \geq u \vee v   \}  \Big )  u^a v^a du dv \, .
\end{multline*} 
So, overall, setting $w_i(x,y):=\p_{x,y}  ( T^* >i, \, W_{i,x}=F_\pi^{-1}
(V_i)   )$,
\begin{multline*}
\p_{\nu \otimes \nu} (T^* > n)  \leq 2 a  \int_0^1 (1-x)^n x^{a-1} dx \\
+ 2 (a+1)^2  \sum_{i=1}^{n-1} \sum_{j=i+1}^{n} \nu\otimes \nu(w_i)  \iint_{[0,1]^2}  \p \Big  (  \bigcap_{k=i +1}^{j} \{  U_k \geq u  \} ,  \bigcap_{k=j +1}^{n} \{  U_k \geq u \vee v   \}  \Big )  u^a v^a du dv \, .
\end{multline*} 
Using the fact that for any $b > -1$,
\begin{equation} \label{majorationevidentegamme}
 \int_0^1 (1-x)^k x^b d x \leq k^{-(b +1)} \int_0^{k} e^{-x} x^b dx \leq k^{-(b +1)} \Gamma (b+1) \, ,
\end{equation}
we get that 
\beq \label{eq1MC}
\p_{\nu \otimes \nu} (T^* > n)   \leq 2 a  \Gamma (a) n^{-a}  
+ 2 (a+1)^2  \sum_{i=1}^{n-1}   \nu\otimes \nu(w_i) \sum_{j=i+1}^{n}    \iint_{[0,1]^2}  (1-u)^{j-i} (1 - u \vee v  )^{n-j} u^a v^a du dv \, .
\eeq
By easy  computations  (that are left to the reader), we infer that Lemma \ref{lmacontinuous} will hold provided one can prove that: 
\begin{Lemma}\label{auxiliary-lem} For any $a>1$, there exists a positive constant $\kappa(a)$ depending only on $a$ such that for any $n \geq 1$,
\beq \label{aimMC2}
\nu\otimes \nu(w_n) \leq \frac{\kappa(a)}{n^a} \, .
\eeq
\end{Lemma}
Obviously, inequality \eqref{aimMC2} holds for any positive integer $n \leq \kappa(a)$.  It is then enough to prove it for $n > \kappa(a)$. Let us do it by recurrence. Hence we assume that for any $k \leq n-1$, 
$\displaystyle \nu\otimes \nu(w_k) \leq \kappa(a) k^{-a}$ and we want to prove it at step $n$. With this aim, we argue as above and infer that 
\[ w_n(x,y) \le (1-y)^n + (a+1)\sum_{i=1}^{n-1}w_i(y,x)\int_{[0,1]}(1-u)^{n-i}u^a du
\, .
\]
Hence, 
\[
\nu\otimes \nu (w_n)\le a \Gamma(a)n^{-a} + (a+1)\sum_{i=1}^{n-1}\nu\otimes \nu (w_i)\int_{[0,1]}(1-u)^{n-i}u^a du \, .
\]
Using the recurrence assumption, it follows that 
\begin{multline*}
\sum_{i=1}^{n-1}\nu\otimes \nu (w_i)\int_{[0,1]}(1-u)^{n-i}u^a du \leq  \sum_{i=1}^{[n/2]}\int_{[0,1]}(1-u)^{n-i}u^a du \\
 +   \frac{ \kappa(a)}{ (  [n/2] +1)^a } \sum_{i=[n/2] +1}^{n-[\log n ]}\int_{[0,1]}(1-u)^{n-i}u^a du \\ +  \frac{ \kappa(a)}{ ( n-[\log n ] +1)^a } \sum_{i=n-[\log n ] +1}^{n -1}\int_{[0,1]}(1-u)^{n-i}u^a du \, .
\end{multline*}
Then, taking into account \eqref{majorationevidentegamme}, we infer that  
\begin{multline*}
(a+1) \sum_{i=1}^{n-1}\nu\otimes \nu (w_i)\int_{[0,1]}(1-u)^{n-i}u^a du  \\ \leq  2^a (a+1) \frac{\Gamma (a)}{ n^{a}} +
 +  2^a  (a+1) \frac{ \kappa(a)}{ n^a [\log n ]^a }  +  a^{-1} \frac{ \kappa(a)}{ ( n-[\log n ] +1)^a }  \, .
\end{multline*}
So, overall, since $n \geq \kappa(a)$, we get 
\[
\nu\otimes \nu (w_n) \le \kappa (a) \rho(a) n^{-a} \, ,
\]
where 
\[
\rho(a) := \big ( a+ 2^a (a+1)  \big ) \frac{a \Gamma(a)}{\kappa (a) } 
 +  \frac{2^a  (a+1) }{  [\log \kappa (a) ]^a }  +  a^{-1}  \Big (  1 -  \frac{\log \kappa (a)}{\kappa (a)}\Big )^{-a} \, .
\]
So choosing $\kappa (a) $ large enough so that $\rho (a) \leq 1$ (which is always possible since $a^{-1} < 1$), inequality \eqref{aimMC2} is proved at step $n$ which ends the recurrence. 
 \hfill $\square$
 
 \subsection{Proof of Lemma \ref{propDR}} \label{sectionproofpropDR}

We start by recalling the inequality line 5 page 27 of 
Dedecker-Rio \cite{DR00}, which holds for every $x,y\in \R$,  
every $n\ge 1$ and any $t >0$: 
\beq \label{ineDeRio}
|W_{n,x}-W_{n,y}|^t \le \alpha^{nt}(|x|+|y|+\Sigma_{n-1})\, |x-y|^t\, ,
\eeq
where $\alpha(u)=1-\frac{C}{(1+u)^{\tau}}$, for every $u\ge 0$,  $\Sigma_0=0$ and 
$\Sigma_{n} =|\varepsilon_1|+\cdots |\varepsilon_n|$, for every $n\ge 1$. 

Denote $\upsilon:=\E(|\varepsilon_1|)$ and let $0<\eta\le 1/\tau-1$.  Notice that $\alpha$ is non-decreasing and bounded by $1$. Hence, for any $n \geq 1$, 
 using that $n\le n^{1/\tau-\eta}$, we get 
\begin{multline*}
\alpha^{nt} (|x|+|y|+\Sigma_{n-1}) \le {\bf 1}_{\{\Sigma_{n-1}>n\upsilon+n^{1/
\tau -\eta}\}}  
+ \alpha^{nt}(2(1+\upsilon)n^{1/\tau -\eta}) {\bf 1}_{\{(1+\upsilon)n^{1/\tau
 -\eta}\ge |x|+|y|\}} \\
+\alpha^{nt}(2(|x|+|y|)){\bf 1}_{\{(1+v) n^{1/\tau -\eta}<|x|+|y|\}}\, .
\end{multline*}
By Theorems 3 and 4 in Baum and Katz \cite{BK63},   since $\mu$ has a moment of order $S$, 
\begin{equation}\label{conv-1}
\sum_{n\ge 1}n^\gamma \P(\Sigma_{n-1}>n\upsilon+n^{1/\tau  -\eta})<\infty\, ,
\end{equation}
provided that $\gamma\le S(1/\tau-\eta)-2$. Since $S/\tau -2\ge 
t/\tau+\gamma$, the latter holds as soon as $\eta\le t/(S\tau)$. 
Hence, we choose $\eta=\min(t/(S\tau), 1/\tau-1)$. On another hand, 
\begin{gather}
\label{conv-2} \sum_{n\ge 1} n^\gamma \alpha^{nt}(2(1+\upsilon)n^{1/\tau -\eta})= \sum_{n\ge 1}
n^\gamma\Big(1-\frac{C}{(1+2(1+\upsilon)n^{1/\tau -\eta})^{\tau}}
\Big)^{nt}<\infty \, .
\end{gather}
Finally,
\begin{gather}
\nonumber\sum_{n\ge 1} n^\gamma\alpha^{nt}(2(|x|+|y|)) {\bf 1}_{\{|x|+|y| >1\}}
= \sum_{n\ge 1}n^\gamma\Big( 1-\frac{C}{ ( 1 + 2
(|x|+|y|) )^{\tau} }
\Big)^{n t } {\bf 1}_{\{|x|+|y| >1\}}\\ 
\label{conv-3}  \le \sum_{n\ge 1}n^\gamma e^{-C n t / (3
(|x|+|y|) )^{\tau}  }  \le D (|x|^{\tau(\gamma+1)}+|y|^{\tau(\gamma+1)})\, ,
\end{gather}
where $D$ is a constant depending on $\gamma$, $t$ and $C$. Starting from \eqref{ineDeRio} and taking into account \eqref{conv-1}, \eqref{conv-2} and \eqref{conv-3} together with the fact that, by \eqref{DR-moment}, 
$\nu$ has a moment of order $S- \tau$ and that $S-\tau \geq \tau(\gamma+1) +t $, we get  the first part of the lemma. 

To prove the last statement, it suffices to notice that for any Lipschitz function $h$ with Lipschitz coefficient equal to $C$, we have, for any $n \geq 2$, 
\[
\delta(n) \leq 2^{-1} C  \sup_{k \geq n-1} \iint \E |W_{k,x} - W_{k,y}| \nu (d x) \nu (dy) \, .
\]
Next simple arguments entail that, for any $n \geq 2$, 
\[
\delta(n)\le C  \iint \E |W_{n-1,x} - W_{n-1,y}| \nu (d x) \nu (dy) \, .
\]
\hfill $\square$

\bigskip

\noindent {\bf Acknowledgement.} The second author is very thankful to 
the laboratories MAP5 and LAMA for their invitations that made possible the present collaboration.

\end{document}